\documentclass[12pt,a4paper,twoside,draft]{article}
\usepackage[english]{babel}
\usepackage[unicode,final,hyperindex]{hyperref}%dvi
\usepackage{cite}
\usepackage{latexsym,amsfonts,amssymb,amsmath,longtable,amsthm}
\usepackage{graphicx,color}
\renewcommand{\le}{\leqslant}
\renewcommand{\ge}{\geqslant}

\begin{document}
\setcounter{page}{3}
\thispagestyle{empty}

\newtheorem*{ttt}{Theorem}
\newtheorem{lem}{Lemma}[section]
\newtheorem{cor}[lem]{Corollary}
\newtheorem{prop}[lem]{Proposition}

\theoremstyle{definition}
\newtheorem{df}[lem]{Definition}
\newtheorem{ex}[lem]{Example}
\newtheorem{tab}{Table}

\theoremstyle{remark}
\newtheorem{rem}[lem]{Remark}
\newtheorem{problem}[lem]{Problem}
\renewcommand{\mod}{\mathrm{mod\ }}

%\hfill УДК 512.542

\begin{center}
%\underline{Sobolev Institute of Mathematics\hspace{4.2cm}Preprint N~225, May
%2009} \vspace{1cm}

{\bf \Large Cocliques of maximal size in the prime graph of a finite simple group}
\vspace{1cm}

{\scshape A. V. Vasil${}^\prime$ev, E. P. Vdovin}
\end{center}

\begin{abstract}
A prime graph of a finite group is defined in the following way: the set of vertices of the graph
is the set of prime divisors of the group, and two distinct vertices $r$ and $s$ are adjacent, if
there is an element of order $rs$ in the group. In this paper we continue our investigation of the
prime graph of a finite simple group started in [1], namely we describe all cocliques of maximal
size for all finite simple groups.
\end{abstract}
\vspace{1cm}

Let $G$ be a finite group, let $\pi(G)$ be the set of all prime divisors of its order, and let
$\omega(G)$ be the spectrum of $G$, i.~e., the set of its element orders. A graph $GK(G)$ is called
the {\em prime graph} (or the {\em Gruenberg-Kegel graph}) of $G$, if the set of vertices of
$GK(G)$ equals $\pi(G)$, and two distinct vertices $r$ and $s$ are adjacent in  if and only if
$rs\in\omega(G)$. Primes $r,s\in \pi(G)$ are called {\em adjacent}, if they are adjacent as
vertices of $GK(G)$. Otherwise, $r$ and $s$ are called {\em non-adjacent}.

In the present paper we continue the investigation of the prime graphs of finite simple groups
started in  \cite{VasVd}. We preserve notation and agreements from~\cite{VasVd}.

Denote by  $t(G)$ the maximal number of prime divisors of $G$ that are pairwise non-adjacent
in~$GK(G)$. In other words $t(G)$ is the maximal number of vertices in cocliques of~$GK(G)$ (a set
of vertices of a graph is called a {\em coclique} (or {\em independent}), if its elements are
pairwise non-adjacent). In the graph theory this number is called an {\em independence number} of a
graph. By analogy we denote by $t(r,G)$ the maximal number of vertices in cocliques of $GK(G)$,
containing a prime~$r$. We call this number  an {\em $r$-independence number}.

In \cite{VasVd} for every finite nonabelian simple group $G$ we gave an arithmetical criterion of
adjacency of vertices in the prime graph $GK(G)$. Using this criterion we determined the values of
$t(G)$, $t(2,G)$, and in case, when $G$ is a group of Lie type over a field of characteristic $p$,
the value of $t(p,G)$. Denote by $\rho(G)$ and $\rho(r,G)$ a coclique of maximal size in $GK(G)$
and a coclique of maximal size, containing $r$, in $GK(G)$ respectively. It is not difficult to see
that in general $\rho(G)$ and $\rho(r,G)$ are not uniquely determined. In \cite{VasVd} all
cocliques $\rho(2,G)$, and also all cocliques $\rho(p,G)$ for groups $G$ of Lie type over a field
of characteristic $p$, were described. Moreover, in the same paper for every simple group $G$ at
least one coclique $\rho(G)$ have been determined, and this allows to calculate $t(G)$, but the
problem of finding all such cocliques has not been considered.

The main goal of the present paper is to find all cocliques of maximal size in the prime graph of a
finite simple group~$G$. In order to achieve this goal we introduce two sets $\Theta(G)$ and
$\Theta'(G)$ consisting of some subsets of $\pi(G)$. Then every coclique $\rho(G)$ of maximal size
can be derived then as $\theta(G)\cup \theta'(G)$, where $\theta(G)\in \Theta(G)$ and
$\theta'(G)\in\Theta'(G)$.

\begin{ttt}\label{MainTheorem} Suppose that $G$ is a finite nonabelian simple group. Then every coclique
of maximal size in $GK(G)$ is the union of $\theta(G)\in\Theta(G)$ and $\theta'(G)\in\Theta'(G)$.
The sets $\Theta(G)$, $\Theta'(G)$ together with the maximal size $t(G)$ of cocliques in $GK(G)$
are described in Proposition~\emph{\ref{indalt}} for alternating groups, in
Table~\emph{\ref{sporadic}} for sporadic groups, and in Tables~{\em\ref{LinearUnitaryTable},
\ref{ClassicTable}, \ref{ExceptTable}} for groups of Lie type.
\end{ttt}

Article \cite{VasVd} appeared to contain several misprints and errors. Some of them were found by
the authors, and others were pointed out by the readers of this paper. We are grateful to W.\,Shi,
H.\,He, A.\,R.\,Moghaddamfar, A.\,Iranmanesh, Z.\,Taheri, S.\,Shariati, M.\,A.\,Grechkoseeva,
A.\,A.\,Buturlakin, and A.\,Zavarnitsine for their comments. Section~\ref{appendix} of the present
article contains the corrections to all detected inaccuracies in~\cite{VasVd}.

\section{Sporadic and alternating groups}\label{sporadicandalt}

Results of the section are easily developed from known ones, and we include them here just
for completeness. Let $G$ be a finite simple sporadic or alternating group. Denote by $\theta(G)$
the intersection of all cocliques of maximal size of $GK(G)$, and by $\Theta(G)$
the set
$\{\theta(G)\}$. The set
$\Theta'(G)$ is defined as follows. A subset $\theta'(G)$ of $\pi(G)\setminus\theta(G)$ is an
element of $\Theta'(G)$ if and only if $\rho(G)=\theta(G)\cup\theta'(G)$ is coclique of $GK(G)$ of
maximal size. Obviously, the sets $\Theta(G)$ and $\Theta'(G)$ are uniquely determined, and
$\Theta'(G)$ either is empty or contains at least two elements.

We start with alternating groups. Let $G=Alt_n$ be the alternating group of degree $n$, and
$n\geqslant5$. In order to describe cocliques of maximal size in $GK(G)$ we present the following
notation. For every prime $r$ define $e(r)=r$ if $r$ is odd, and $e(r)=4$ if $r=2$. Denote by
$\tau(n)$ the set of all primes $r$ with $n/2\leqslant e(r)\leqslant n$, and by $s_n$ and $s'_n$
the smallest elements of~$\tau(n)$ and $\tau(n)\setminus\{s_n\}$ respectively. Define the sets
$\tau'(n)$ and $\tau''(n)$ as follows. A prime $r$ lies in $\tau'(n)$ if and only if $e(r)<n/2$ and
$e(r)+e(s_n)>n$, and $r$ lies in $\tau''(n)$ if and only if $e(r)<n/2$ and $e(r)+e(s'_n)>n$.

\begin{prop}\label{indalt}
Let $G$ be an alternating group of degree $n$, and $n\geqslant5$.
\begin{itemize}
\item[{\em 1.}] If $\tau'(n)=\tau''(n)=\varnothing$, then $\theta(G)=\tau(n)$, and $\Theta'(G)=\varnothing$.

\item[{\em 2.}] If $\tau'(n)=\varnothing$ and $\tau''(n)\neq\varnothing$, then
$\theta(G)=\tau(n)\setminus\{s_n\}$, and $\Theta'(G)=\{\{r\}\mid r\in\tau''(n)\cup\{s_n\}\}$.

\item[{\em 3.}] If $|\tau'(n)|=1$, then $\theta(G)=\tau(n)\cup\tau'(n)$, and $\Theta'(G)=\varnothing$.

\item[{\em 4.}] If $|\tau'(n)|\geqslant2$, then $\theta(G)=\tau(n)$, and $\Theta'(G)=\{\{r\}\mid r\in \tau'(n)\}$.
\end{itemize}

In all cases every coclique of maximal size in $GK(G)$ is of the form $\theta(G)\cup\theta'(G)$,
where $\theta(G)\in\Theta(G)$, and $\theta'(G)\in\Theta'(G)$. The set $\Theta(G)=\{\theta(G)\}$ is
one-element, and all elements $\theta'(G)$ of $\Theta'(G)$ are one-element subsets of $\pi(G)$.
\end{prop}

\begin{proof}
An adjacency criterion for vertices of $GK(G)$ \cite[Proposition~1.1]{VasVd} can be formulated as
follows. Distinct primes $r,s\in\pi(G)$ are non-adjacent in $GK(G)$ if and only if $e(r)+e(s)>n$.
Therefore, $\tau(n)$ is a coclique of $GK(G)$, and $\pi(G)\setminus\tau(n)$ is a clique. Moreover,
if $r,s,t$ are distinct primes from $\pi(G)$, and $e(r)<e(s)$, then the adjacency of $s$ and $t$
implies the adjacency of $r$ and $t$, as well as the non-adjacency of $r$ and $t$ implies the
non-adjacency of $s$ and $t$. These simple observations allow to verify the assertion easily.
\end{proof}

\begin{prop}\label{indsporadic}
Let $G$ be a simple sporadic group. If $\Theta'(G)=\varnothing$, then  $\theta(G)$ is the unique
coclique of maximal size in $GK(G)$. If $\Theta'(G)\neq\varnothing$, then every coclique of maximal
size is of the form $\theta(G)\cup\theta'(G)$, where $\theta'(G)\in\Theta'(G)$. If $G\neq M_{23}$
then every $\theta'(G)$ of $\Theta'(G)$ contains precisely one element. The sets $\Theta(G)$ and
$\Theta'(G)$, as well as the value of $t(G)$, are listed in
Table~{\em\ref{sporadic}}.
\end{prop}

\begin{proof}
The proposition is easy to verify using
\cite{Atlas} or~\cite{GAP}.
\end{proof}

\textsl{Remark.} Note that in Columns 3 and 4 of Table 1 we list the elements of $\Theta(G)$ and
$\Theta'(G)$, that is sets $\theta(G)\in\Theta(G)$ and $\theta'(G)\in\Theta'(G)$, and omit the
braces for one-element sets. In particular, for group $G=M_{11}$ we have
$\Theta(G)=\{\theta(G)\}=\{\{5,11\}\}$ and $\Theta'(G)=\{\{2\},\{3\}\}$, while for $G=M_{23}$ we
have $\Theta(G)=\{\theta(G)\}=\{\{11,23\}\}$ and $\Theta'(G)=\{\{2,5\},\{3,7\}\}$.

\begin{tab}\label{sporadic}{\bfseries Cocliques of sporadic groups}\vspace{1\baselineskip}

{\small
\begin{tabular}{|r|c|l|l|}
 \hline $G$ & $t(G)$ & $\Theta(G)$ & $\Theta'(G)$ \\
 \hline $M_{11}$ & $3$ & $\{5,11\}$ & $2,3$  \\
 $M_{12}$ & $3$ & $\{3,5,11\}$ & $\varnothing$   \\
 $M_{22}$ & $4$ & $\{5,7,11\}$ & $2,3$  \\
 $M_{23}$ & $4$ & $\{11,23\}$ & $\{2,5\}$, $\{3,7\}$   \\
 $M_{24}$ & $4$ & $\{5,7,11,23\}$ & $\varnothing$   \\
 $J_{1}$ & $4$ & $\{7,11,19\}$ & $2,3,5$  \\
 $J_{2}$ & $2$ & $7$ & $2,3,5$  \\
 $J_{3}$ & $3$ & $\{17,19\}$ & $2,3,5$   \\
 $J_{4}$ & $7$ & $\{11,23,29,31,37,43\}$ & $5,7$  \\
 $\operatorname{Ru}$ & $4$ & $\{7,13,29\}$ & $3,5$  \\
 $\operatorname{He}$ & $3$ & $\{5,7,17\}$ & $\varnothing$  \\
 $\operatorname{McL}$ & $3$ & $\{7,11\}$ & $3,5$ \\
 $\operatorname{HN}$ & $3$ & $\{11,19\}$ & $3,5,7$ \\
 $\operatorname{HiS}$ & $3$ & $\{7,11\}$ & $2,3,5$ \\
 $\operatorname{Suz}$ & $4$ & $\{5,7,11,13\}$ & $\varnothing$ \\
 $\operatorname{Co}_{1}$ & $4$ & $\{11,13,23\}$ & $5,7$ \\
 $\operatorname{Co}_{2}$ & $4$ & $\{7,11,23\}$ & $3,5$ \\
 $\operatorname{Co}_{3}$ & $4$ & $\{5,7,11,23\}$ & $\varnothing$ \\
 $\operatorname{Fi}_{22}$ & $4$ & $\{5,7,11,13\}$ & $\varnothing$ \\
 $\operatorname{Fi}_{23}$ & $5$ & $\{11,13,17,23\}$ & $5,7$ \\
 $\operatorname{Fi}'_{24}$ & $6$ & $\{11,13,17,23,29\}$ & $5,7$ \\
 $\operatorname{O'N}$ & $5$ & $\{7,11,19,31\}$ & $3,5$ \\
 $\operatorname{LyS}$ & $6$ & $\{5,7,11,31,37,67\}$ & $\varnothing$ \\
 $F_{1}$ & $11$ & $\{11,13,19,23,29,31,41,47,59,71\}$ & $7,17$ \\
 $F_{2}$ & $8$ & $\{7,11,13,17,19,23,31,47\}$ & $\varnothing$ \\
 $F_{3}$ & $5$ & $\{5,7,13,19,31\}$ & $\varnothing$ \\
   \hline
\end{tabular}}
\end{tab}

In addition, we notice another substantial property of prime graphs of groups under consideration.

\begin{prop}\label{cliquesporalt} Suppose that $G$ is either an alternating group of degree $n$, $n\geqslant5$,
or a sporadic group distinct from $M_{23}$. Then the set $\pi(G)\setminus\theta(G)$ is a clique of
$GK(G)$.
\end{prop}

\begin{proof}
This follows from \cite{Atlas} and \cite[Proposition~1.1]{VasVd}.
\end{proof}

\section{Preliminary results for groups of Lie type}\label{preliminary}

We write $[x]$ for the integer part of a rational number $x$. The set of prime divisors of a
natural number $m$ is denoted by $\pi(m)$. By  $(m_1,m_2,\dots,m_s)$ we denote the  greatest common
divisor of numbers $m_1,m_2,\dots,m_s$. For a natural number $r$, the $r$-part of a natural number
$m$ is the greatest divisor $t$ of $m$ with $\pi(t)\subseteq\pi(r)$. We write $m_r$ for the
$r$-part of $m$ and  $m_{r'}$ for the quotient $m/m_r$.

If $q$ is a natural number, $r$ is an odd prime and $(q,r)=1$, then $e(r,q)$ denotes
a~multi\-plicative order of $q$ modulo $r$, that is a minimal natural number $m$ with
$q^m\equiv1\pmod{r}$. For an odd $q$, we put $e(2,q)=1$ if $q\equiv1\pmod{4}$, and $e(2,q)=2$
otherwise.

\begin{lem}\label{Zsigmondy Theorem}
{\em (Corollary of Zsigmondy's theorem \cite{zs})} Let $q$ be a natural number greater than $1$.
For every natural number $m$ there exists a prime $r$ with $e(r,q)=m$ but for the cases $q=2$ and
$m=1$, $q=3$ and $m=1$, and $q=2$ and $m=6$.
\end{lem}

\textsl{Remark.} In conclusion of the same corollary \cite[Lemma~1.4]{VasVd} in our previous
article we miss two exceptions: $m=1$ and $q=2$, and $m=1$ and $q=3$. However, these exceptions
don't arise in all proofs and arguments from \cite{VasVd}, that use the corollary to Zsigmondy's
theorem.

\smallskip

A prime $r$ with $e(r,q)=m$ is called a {\em primitive prime divisor} of $q^m-1$. By Lemma
\ref{Zsigmondy Theorem} such a number exists except for the cases mentioned in the lemma. Given $q$
we denote by $R_m(q)$ the set of all primitive prime divisors of $q^m-1$ and by $r_m(q)$ any
element of $R_m(q)$. If $m\not=2$ then a divisor $k_m(q)$ of $q^m-1$ is said to be the {\em
greatest primitive divisor} if $\pi(k_m(q))=\pi(R_m(q))$ and $k_m(q)$ is the greatest divisor with
this property, i.e., $k_m(q)=(q^m-1)_t$, where $t=\prod_{s\in R_m(q)}s$. The {\em greatest
primitive divisor} $k_2(q)$ of $q^2-1$ is the greatest divisor of $q+1$ with $\pi(k_2(q))=R_2(q)$.
The singularity in the definition of the greatest primitive divisor in case $m=2$ appears because
of the singularity of the definition for $e(2,q)$. Following our definition of $e(2,q)$, we derive
that $k_1(q)=(q-1)/2$ if $q\equiv-1\pmod{4}$, and $k_1(q)=q-1$ otherwise; $k_2(q)=(q+1)/2$ if
$q\equiv1\pmod{4}$, and $k_2(q)=q+1$ otherwise. The following lemma provides a formula for
expressing greatest primitive divisors $k_m$, $m\geqslant3$, in terms of cyclotomic polynomials
$\phi_m(x)$.

\begin{lem}\label{gpd} {\em \cite{R}}
Let $q$ and $m$ be natural numbers, $q>1$, $m\geqslant3$, and let $k_m(q)$ be the greatest primitive
divisor of $q^m-1$. Then $$k_m(q)=\frac{\phi_m(q)}{( \phi_{m_{r'}}(q),r)},$$
where $r$ is the greatest prime divisor of~$m$.
\end{lem}

Usually the number $q$ is fixed (for example, by the
choice of a group of Lie type $G$), and we write $R_m$, $r_m$, and $k_m$ instead
of $R_m(q)$, $r_m(q)$, and $k_m(q)$ respectively.
According to our definitions, if $i\neq j$, then
$\pi(R_i)\cap\pi(R_j)=\varnothing$, and so
$(k_i,k_j)=1$.

%\smallskip

\begin{lem}\label{Divisibility} {\em \cite[Lemma~6(iii)]{ZavL3}} Let $q,k,l$ be natural numbers. Then

\emph{(a)} $(q^k-1,q^l-1)=q^{(k,l)}-1;$

\emph{(b)} $(q^k+1,q^l+1)=\left \{
\begin{array}{ll}
q^{(k,l)}+1, &\mbox{if both $\frac{k}{(k,l)}$ and $\frac{l}{(k,l)}$ are
odd},\\
(2,q+1), & \mbox{otherwise};
\end{array} \right.$

\emph{(c)} $(q^k-1,q^l+1)=\left \{
\begin{array}{ll}
q^{(k,l)}+1,&\mbox{if $\frac{k}{(k,l)}$ is even and $\frac{l}{(k,l)}$ is
odd},\\
(2,q+1), & \mbox{otherwise.}
\end{array} \right.$

In particular, for every $q\ge2$, $k\ge 1$ the inequality $(q^k-1,q^k+1)\le 2$ holds.
\end{lem}

%\begin{cor}\label{cornum} Let $q,k,l$ be natural numbers, and $q>1$. Then

%\emph{(a)} $q^k-1$ divides $q^l-1$ iff $k$ divides $l;$

%\emph{(b)} $q^k+1$ divides $q^l+1$ iff $l/k$ is an odd natural number$;$

%\emph{(c)} $q^k+1$ divides $q^l-1$ iff $l/k$ is an even natural number$;$

%\emph{(d)} $q^k-1$ divides $q^l+1$ iff $q=2,3$ and $k=1$, or $q=2$, $k=2$ and $l$ is odd.

%\end{cor}

We recall also the following statements  \cite[statements~(1) and~(4)]{VasVd}.
Given $q=p^\alpha$, where $p$ is a prime, and odd prime $c\not=p$ we have:

\begin{equation}\label{firstold}
c\text{ divides }q^x-1\text{ if and only if }\\
e(c,q)\text{ divides }x;
\end{equation}

\begin{equation}\label{fourthold}
\text{if }c\text{ divides }q^x-\epsilon \text{, where
}\epsilon\in\{+1,-1\},\\ \text{ then } \eta(e(c,q))\text{
divides }x.
\end{equation}

The function $\eta(n)$  is defined in Proposition~\ref{adjbn}.

In the proofs of Propositions \ref{adjbn}, \ref{adjdn}, and \ref{adjexcept} by $\epsilon,
\epsilon_i$ we denote elements from the set $\{+1,\-1\}$. For groups of Lie type our notation
agrees with that of \cite{VasVd}. We write $A_n^\varepsilon(q)$, $D_n^\varepsilon(q)$, and
$E_6^\varepsilon(q)$, where $\varepsilon\in\{+,-\}$, and $A_n^+(q)=A_n(q)$, $A_n^-(q)={}^2A_n(q)$,
$D_n^+(q)=D_n(q)$, $D_n^-(q)={}^2D_n(q)$, $E_6^+(q)=E_6(q)$, $E_6^-(q)={}^2E_6(q)$. In
\cite[Proposition~2.2]{VasVd}, considering unitary groups, we define the function
\begin{equation}\label{nu(m)}
\nu(m)=\left\{
\begin{array}{rl}
m &\text{ if }m\equiv 0(\mod 4),\\
\frac{m}{2}& \text{ if }m\equiv 2(\mod 4),\\
2m&\text{ if }m\equiv1(\mod 2).\\
\end{array}\right.
\end{equation}
Clearly $\nu(m)$ is a bijection from $\mathbb{N}$ onto $\mathbb{N}$ and $\nu^{-1}(m)=\nu(m)$. In most cases it is natural to consider linear
and unitary groups together. So we define
\begin{equation}\label{nuepsilon(n)}
\nu_{\varepsilon}(m)=\left\{
\begin{array}{rl}
m &\text{ if } \varepsilon=+,\\
\nu(m) & \text{ if }\varepsilon=-.\\
\end{array}\right.
\end{equation}

\begin{prop}\label{adjbn}
Let $G$ be one of simple groups of Lie type $B_n(q)$ or $C_n(q)$  over a field of
characteristic~$p$. Define
$$\eta(m)=\left\{
\begin{array}{cc}
m &\text{ if }m\text{ is odd},\\
\frac{m}{2}& \text{ otherwise}.\\
\end{array}\right.$$ Let $r,s$ be odd primes with $r,s\in\pi(G)\setminus\{p\}$. Put $k=e(r,q)$ and $l=e(s,q)$, and suppose that $1\le
\eta(k)\le \eta(l)$. Then $r$ and $s$ are non-adjacent if and only if  $\eta(k)+\eta(l)> n$, and
$k$, $l$ satisfy the following condition{\em:}
\end{prop}

\begin{equation}\label{strange}
\dfrac{l}{k} \text{ {\em is not an odd natural number}}
\end{equation}

\begin{proof}
We prove the ``if'' part first. Assume that  $\eta(k)+\eta(l)\le n$, then
there exists a maximal torus $T$ of order $\frac{1}{(2,q-1)}(q^{\eta(k)}+(-1)^k)(q^{\eta(l)}+(-1)^l)(q-1)^{n-\eta(k)-\eta(l)}$  of $G$ (see
\cite[Lemma~1.2(2)]{VasVd}, for example).
Both $r,s$ divide $\vert T\vert$, hence $r,s$ are adjacent in~$G$. If
$\frac{l}{k}$ is an odd integer, then either both $k,l$ are odd and Lemma~\ref{Divisibility}(a) implies that
$q^{\eta(k)}+(-1)^k=q^k-1$ divides $q^{\eta(l)}+(-1)^l=q^l-1$, or both $k,l$ are even and Lemma~\ref{Divisibility}(b) implies that
$q^{\eta(k)}+(-1)^k=q^{k/2}+1$ divides $q^{\eta(l)}+(-1)^l=q^{l/2}+1$. Again both $r,s$  divide $\vert T\vert$, where $T$ is a maximal
torus of order $\frac{1}{(2,q-1)}(q^{\eta(l)}+(-1)^l)(q-1)^{n-\eta(l)}$ of $G$ (the existence of such torus
follows from \cite[Lemma~1.2(2)]{VasVd}), so $r,s$ are adjacent.

Now we prove the ``only if'' part. Assume by contradiction that
$\eta(k)+\eta(l)>n$ and $l/k$ is not an odd natural
number, but $r,s$
are adjacent. Then $G$ contains an element $g$ of order $rs$. The element $g$ is semisimple, since $(rs,p)=1$, hence $g$ is contained in a
maximal torus $T$ of $G$. By \cite[Lemma~1.2(2)]{VasVd} it follows that $\vert
T\vert=\frac{1}{(2,q-1)}(q^{n_1}-\epsilon_1)(q^{n_2}-\epsilon_2)\ldots(q^{n_k}-\epsilon_k)$, where $n_1+n_2+\ldots+n_k=n$.
Up to renumberring, we may assume
that $r$ divides $(q^{n_1}-\epsilon_1)$, while $s$ divides either $(q^{n_1}-\epsilon_1)$ or~$(q^{n_2}-\epsilon_2)$. Assume first that $s$
divides~$(q^{n_2}-\epsilon_2)$. Then \eqref{fourthold} implies that $\eta(k)$ divides $n_1$ and $\eta(l)$ divides $n_2$, so $n_1+n_2\ge
\eta(k)+\eta(l)>n$, a contradiction.

Now assume that both $r,s$ divide $(q^{n_1}-\epsilon_1)$. Again \eqref{fourthold} implies that both $\eta(k),\eta(l)$ divide $n_1$. Now
$\eta(k)+\eta(l)>n$ and $\eta(k)\le \eta(l)$, so $\eta(l)=n_1$. Assume first that $l$ is odd. Then $l=\eta(l)=n_1$ and $s$ divides $q^l-1$.
Since $s$ is odd, Lemma \ref{Divisibility} imples that $s$ does not divide $q^l+1$, hence $q^{n_1}-\epsilon_1=q^{n_1}-1$. Since $r$ divides
$q^{n_1}-1$, by using \eqref{firstold} we obtain that $k$ divides $n_1=l$, hence $k$ is odd. Therefore
$\frac{l}{k}$ is an odd integer, a contradiction with~\eqref{strange}. Now assume that $l$ is even. Then $l/2=\eta(l)=n_1$
and $s$ divides $q^l-1$. In view of \eqref{firstold}, $s$ does not divide $q^{l/2}-1$, hence $s$ divides $q^{l/2}+1$ and
$q^{n_1}-\epsilon_1=q^{n_1}+1$. Now \eqref{fourthold} implies that $\eta(k)$ divides $n_1$, hence $k$ divides $2n_1=l$. By Lemma
\ref{Divisibility}(c) we obtain that $r$ does not divide $q^{l/2}-1$, hence $k$ does not divide $l/2$ and $\frac{l}{k}$ is an
odd integer, a contradiction with~\eqref{strange}.
\end{proof}

\begin{prop}\label{adjdn}
Let $G=D_n^{\varepsilon}(q)$ be a finite simple group of Lie type over a field of characteristic  $p$, and let the function
$\eta(m)$ be defined as in Proposition~{\em\ref{adjbn}}. Suppose $r,s$ are odd primes and $r,s\in\pi(D_n^\varepsilon(q))\setminus\{p\}$.
Put $k=e(r,q)$, $l=e(s,q)$, and  $1\le\eta(k)\le\eta(l)$. Then $r$ and $s$ are non-adjacent if and only if $2\cdot\eta(k)+2\cdot\eta(l)>
2n-(1-\varepsilon(-1)^{k+l})$,
$k$ and $l$ satisfy~\eqref{strange}, and, if  $\varepsilon=+$, then the chain of equalities{\em:}
\begin{equation}\label{strange2}
n=l=2\eta(l)=2\eta(k)=2k
\end{equation}
is not true.
\end{prop}

\begin{proof}
The following inclusions are known $\widetilde{B}_{n-1}(q)\leq \widetilde{D}_n^\varepsilon(q)\leq
\widetilde{B}_n(q)$ (see \cite[Table~2]{KondSbgrps}), where $\widetilde{B}_{n-1}(q)$,
$\widetilde{D}_n^\varepsilon(q)$, $\widetilde{B}_n(q)$ are central extensions of corresponding
simple groups and $n\ge 4$. Since the Schur multiplier for each of simple groups $B_{n-1}(q)$,
$D_n^\varepsilon(q)$, $B_n(q)$ has order equal to $1$, $2$, or $4$,  it is clear that two odd
prime divisors of the order of a simple group isomorphic to $B_n(q)$ or $D_n^\varepsilon(q)$ are
adjacent if and only if they are adjacent in every central extension of the group. Hence if two
odd prime divisors of $\vert D_n^\varepsilon(q)\vert$ are adjacent in $GK(B_{n-1}(q))$, then they
are adjacent in $GK(D_n^\varepsilon(q))$ and if two odd prime divisors of $\vert
D_n^\varepsilon(q)\vert$ are non-adjacent in $GK(B_{n}(q))$, then they are non-adjacent in
$GK(D_n^\varepsilon(q))$. There can be the following cases:
\begin{itemize}
\item[(i)] $\eta(k)+\eta(l)\le n-1$;
\item[(ii)]$\eta(k)+\eta(l)\ge n$, $l/k$ is an odd number and $\eta(l)\le n-1$;
\item[(iii)] $\eta(k)+\eta(l)=n$ and $\frac{l}{k}$ is not an odd natural number;
\item[(iv)]  $\eta(l)=n$ and $\frac{l}{k}$ is an odd natural number;
\item[(v)] $\eta(k)+\eta(l)>n$ and  $\frac{l}{k}$ is not an odd natural number.
\end{itemize}
By Lemma \ref{adjbn} in cases (i), (ii) primes $r,s$ are adjacent in $GK(B_{n-1}(q))$, while in case (v) primes $r,s$ are non-adjacent in
$GK(B_n(q))$. In view of above notes it follows that we
need to consider (iii) and~(iv).

Assume first that $\eta(k)+\eta(l)=n$ and $\frac{l}{k}$ is not an odd natural
number, i.~e., case (iii) holds. Since $(rs,p)=1$, the primes  $r,s$ are
adjacent in
$GK(D_n^\varepsilon(q))$ if and only if there exists  a maximal torus $T$ of
$G$ of order divisible by $rs$.  In view of \cite[Lemma~1.2(3)]{VasVd} the order
$\vert T\vert$ is equal to
$\frac{1}{(4,q^n-\varepsilon1)}(q^{n_1}-\epsilon_1)\cdot\ldots\cdot(q^{n_m}
-\epsilon_m)$, where $n_1+\ldots+n_m=n$  and
$\epsilon_1\cdot\ldots\cdot\epsilon_m=\varepsilon1$.
Up to renumberring, we may assume that $r$ divides
$q^{n_1}-\epsilon_1$, while $s$ divides either $q^{n_1}-\epsilon_1$, or $q^{n_2}-\epsilon_2$.

If $s$ divides $q^{n_1}-\epsilon_1$, then
\eqref{fourthold}
implies that both $\eta(k)$, $\eta(l)$ divide $n_1$. As in the proof of Proposition \ref{adjbn} we derive that $r,s$ are adjacent if and
only if $\frac{l}{k}$ is an
odd integer.

Assume now that $s$ divides $q^{n_2}-\epsilon_2$. Then \eqref{fourthold} implies that $\eta(k)$ divides $n_1$ and $\eta(l)$ divides $n_2$.
Hence we obtain
the following inequalities $n\ge n_1+n_2\ge \eta(k)+\eta(l)=n$, so $\eta(k)=n_1$, $\eta(l)=n_2$, and $q^{n_1}-\epsilon_1=q^{\eta(k)}+(-1)^k$,
$q^{n_2}-\epsilon_2=q^{\eta(l)}+(-1)^l$. If $\varepsilon=-$, then a maximal torus $T$ of order
$\frac{1}{(4,q^n+1)}(q^{\eta(k)}+(-1)^k)(q^{\eta(l)}+(-1)^l)$
of $G$ exists if and only if $k,l$ have the distinct parity, i.~e., if and only
if $2n-(1-\varepsilon(-1)^{k+l})=2n-(1+(-1)^{k+l})=2n$. Hence
in
this case
$r,s$ are non-adjacent if and only if the inequality $2\cdot\eta(k)+2\cdot\eta(l)> 2n-(1-\varepsilon(-1)^{k+l})$ holds. If $\varepsilon=+$ and
$n_1\not=n_2$, then a maximal torus $T$ of order
$\frac{1}{(4,q^n-1)}(q^{\eta(k)}+(-1)^k)(q^{\eta(l)}+(-1)^l)$ of $G$ exists if
and only if $k,l$ have
the same parity, i.~e., if and only if $2n-(1-\varepsilon(-1)^{k+l})=2n-(1-(-1)^{k+l})=2n$.  Hence in this case $r,s$ are non-adjacent if and only if
the inequality $2\cdot\eta(k)+2\cdot\eta(l)> 2n-(1-\varepsilon(-1)^{k+l})$ holds. If $n_1=n_2=n/2$ and $\frac{l}{k}$ is an odd integer,
then, $r,s$ are adjacent. Assume that $n_1=n_2=n/2$ and $\frac{l}{k}$ is not an odd integer. The condition $\frac{l}{k}$ is not an odd
integer implies that $l\not=k$, so the chain of equalities \eqref{strange2} holds.
In this case there exists a maximal torus $T$ of order
$\frac{1}{(4,q^n-1)}(q^n-1)=\frac{1}{(4,q^n-1)}(q^{n/2}-1)(q^{n/2}+1)$ of $G$, so condition \eqref{strange2} is not satisfied and $r,s$ are adjacent.

Now assume that $\eta(l)=n$ and $\frac{l}{k}$ is an odd natural number, i.~e., case (iv) holds. In this case there exists a maximal torus
$T$ of order
$\frac{1}{(4,q^n-\varepsilon1)}(q^n+(-1)^l)$ of $G$ (if such a torus does not exist then $s$ does not divide $\vert G\vert$). The fact that
$\frac{l}{k}$ is an odd prime implies that $r$ divides $\vert T\vert$, so $r,s$ are adjacent.
\end{proof}

Now we consider simple exceptional groups of Lie type. Note that the orders of maximal tori of
simple exceptional groups were listed in \cite[Lemma~1.3]{VasVd}. However, for groups $E_7(q)$,
$E_8(q)$, and Ree groups ${}^2F_4(2^{2n+1})$ (items (4), (5), and (9) of the lemma respectively),
the list of orders of tori is incorrect. The following lemma corrects this.

\begin{lem}\label{toriofexcptgrps} {\em (see \cite{Ca3}% and \cite{der1}
)} Let $\overline{G}$ be a connected simple exceptional algebraic
group
of adjoint type and let $G=O^{p'}(\overline{G}_\sigma)$ be a finite simple
exceptional group of
Lie type.

\item[{\em 1.}] For every maximal torus $T$ of $G=E_7(q)$, the number
  $m=(2,q-1)\vert T\vert$ is equal to one of the following:
$(q+1)^{n_1}(q-1)^{n_2},$ $n_1+n_2=7;$ $(q^2+1)^{n_1}(q+1)^{n_2}(q-1)^{n_3},$ $1\leqslant
n_1\leqslant2,$ $2n_1+n_2+n_3=7,$ and $m\neq(q^2+1)(q\pm1)^5;$
$(q^3+1)^{n_1}(q^3-1)^{n_2}(q^2+1)^{n_3}(q+1)^{n_4}(q-1)^{n_5},$ $1\leqslant n_1+n_2\leqslant2,$
$3n_1+3n_2+2n_3+n_4+n_5=7,$ and $m\neq(q^3+\epsilon1)(q-\epsilon1)^4,$ $m\neq(q^3\pm1)(q^2+1)^2,$
$m\neq(q^3+\epsilon1)(q^2+1)(q+\epsilon1)^2;$
$(q^4+1)(q^2\pm1)(q\pm1);$ $(q^5\pm1)(q^2-1);$ $(q^5+\epsilon1)(q+\epsilon1)^2;$ $q^7\pm1;$
$(q-\epsilon1)\cdot (q^2+\epsilon q+1)^3; (q^5-\epsilon1)\cdot
  (q^2+\epsilon q +1); (q^3\pm1)\cdot (q^4-q^2+1); (q-\epsilon1)\cdot
  (q^6+\epsilon q^3+1);$ $(q^3-\epsilon1)\cdot(q^2-\epsilon
  q+1)^2,$ where $\epsilon=\pm$. Moreover, for every number $m$ given above there exists a torus $T$ with $(2,q-1)\vert T\vert=m$.

\item[{\em 2.}] Every maximal torus $T$ of $G=E_8(q)$ has one of the
  following orders:
$(q+1)^{n_1}(q-1)^{n_2},$ $n_1+n_2=8;$ $(q^2+1)^{n_1}(q+1)^{n_2}(q-1)^{n_3},$ $1\leqslant
n_1\leqslant4,$ $2n_1+n_2+n_3=8,$ and $|T|\neq(q^2+1)^3(q\pm1)^2,$ $|T|\neq(q^2+1)(q\pm1)^6;$
$(q^3+1)^{n_1}(q^3-1)^{n_2}(q^2+1)^{n_3}(q+1)^{n_4}(q-1)^{n_5},$ $1\leqslant n_1+n_2\leqslant2,$
$3n_1+3n_2+2n_3+n_4+n_5=8,$ and $|T|\neq(q^3\pm1)^2(q^2+1),$
$|T|\neq(q^3+\epsilon1)(q-\epsilon1)^5,$ $|T|\neq(q^3+\epsilon1)(q^2+1)(q+\epsilon1)^3,$
$|T|\neq(q^3+\epsilon1)(q^2+1)^2(q-\epsilon1);$ $q^8-1;$ $(q^4+1)^2;$ $(q^4+1)(q^2\pm1)(q\pm1)^2;$
$(q^4+1)(q^2-1)^2;$ $(q^4+1)(q^3+\epsilon1)(q-\epsilon1);$
$(q^5+\epsilon1)(q+\epsilon1)^3;$ $(q^5\pm1)(q+\epsilon1)^2(q-\epsilon1);$
$(q^5+\epsilon1)(q^2+1)(q-\epsilon1);$ $(q^5+\epsilon1)(q^3+\epsilon1);$ $(q^6+1)(q^2\pm1);$
$(q^7\pm1)(q\pm1);$ $(q-\epsilon1)\cdot (q^2+\epsilon q+1)^3\cdot(q\pm1);$ $(q^5-\epsilon1)\cdot
  (q^2+\epsilon q +1)\cdot(q+\epsilon1);$ $(q^3\pm1)\cdot (q^4-q^2+1)\cdot(q\pm1);$
$(q-\epsilon1)\cdot
  (q^6+\epsilon q^3+1)\cdot(q\pm1);$ $(q^3-\epsilon1)\cdot(q^2-\epsilon
  q+1)^2\cdot(q\pm1);$ $q^8-q^4+1;$
  $q^8+q^7-q^5-q^4-q^3+q+1;$
 $q^8-q^6+q^4-q^2+1;$ $(q^4-q^2+1)^2;$ $(q^6+\epsilon q^3+1)(q^2+\epsilon q+1);$
 $q^8-q^7+q^5-q^4+q^3-q+1;$ $(q^4+\epsilon q^3+q^2+\epsilon q+1)^2;$
  $(q^4-q^2+1)(q^2\pm q+1)^2;$
   $(q^2-q+1)^2\cdot(q^2+q+1)^2;$
 $(q^2\pm q+1)^4,$ where $\epsilon=\pm$. Moreover, for every number
given above there exists a torus of corresponding order.
\item[{\em 3.}] Every maximal torus $T$ of $G={^2F_4(2^{2n+1})}$ with $n\ge1$
has one of the following orders: $q^2+\epsilon q\sqrt{2q}+q+\epsilon \sqrt{2q}+1;$ $q^2-\epsilon
q\sqrt{2q}+\epsilon \sqrt{2q}-1;$ $q^2-q+1;$ $(q\pm \sqrt{2q}+1)^2;$ $(q-1)(q\pm \sqrt{2q}+1);$
$(q\pm1)^2;q^2\pm1;$ where $q=2^{2n+1}$ and~$\epsilon=\pm$. Moreover, for every number given above
there exists a torus of corresponding order.
\end{lem}

\begin{prop}\label{adjexcept}
Let $G$ be a finite simple exceptional group of Lie type over a field of order
$q$ and  characteristic~$p$. Suppose that
$r,s$ are odd primes, and assume that  $r,s\in\pi(G)\setminus\{p\}$. Put
$k=e(r,q)$, $l=e(s,q)$, and assume that
$1\le k\le l$. Then  $r$ and $s$ are non-adjacent if and only if  $k\not=l$ and one of the following holds{\em:}
\begin{itemize}
\item[{\em 1.}] $G=G_2(q)$  and either $r\not=3$ and $l\in\{3,6\}$ or $r=3$
and~${l=9-3k}$.
\item[{\em 2.}] $G=F_4(q)$ and either $l\in\{8,12\}$, or $l=6$
and~$k\in\{3,4\}$, or $l=4$ and~${k=3}$.
\item[{\em 3.}] $G=E_6(q)$ and either $l=4$ and $k=3$, or $l=5$ and $k\ge3$, or $l=6$
and $k=5$, or $l=8$, $k\ge3$, or $l=8$, $r=3$, and $(q-1)_3=3$, or
$l=9$, or $l=12$ and $k\not=3$.
\item[{\em 4.}] $G={^2E_6(q)}$ and either $l=6$ and $k=4$, or $l=8$, $k\ge3$, or
$l=8$, $r=3$, and
$(q+1)_3=3$, or  $l=10$ and $k\ge 3$, or $l=12$ and $k\not=6$, or~$l=18$.
\item[{\em 5.}] $G=E_7(q)$ and either $l=5$ and $k=4$, or  $l=6$ and $k=5$, or $l\in\{14,18\}$ and
$k\not=2$, or  $l\in\{7,9\}$ and $k\ge2$, or $l=8$ and $k\ge3,k\not=4$,  or $l=10$ and $k\ge3,
  k\not=6$, or $l=12$ and~$k\ge 4,k\not=6$.
\item[{\em 6.}] $G=E_8(q)$ and either $l=6$ and $k=5$, or $l\in\{7,14\}$ and $k\ge3$, or
$l=9$ and $k\ge 4$, or $l\in\{8,12\}$ and $k\ge 5, k\not=6$, or $l=10$ and $k\ge3, k\not=4,6$, or $l=18$ and $k\not=1,2,6$, or
$l=20$ and $r\cdot k\not=20$, or~$l\in\{15,24,30\}$.
\item[{\em 7.}] $G={^3D_4(q)}$ and either $l=6$ and $k=3$,
or~$l=12$.
\end{itemize}
\end{prop}

\begin{proof}
Recal that $k_m$ is the greatest primitive divisor of $q^m-1$, while $R_m$ is the set of all prime primitive divisors of $q^m-1$.
The orders of maximal tori in exceptional groups are given
in~\cite[Lemma~1.3]{VasVd} and Lemma \ref{toriofexcptgrps}, for example.

1. Since $\vert G_2(q)\vert=q^6(q^2-1)(q^6-1)$, the numbers $k,l$ are in the set $\{1,2,3,6\}$.  If $\{k,l\}\subseteq \{1,2\}$, then the
existence
of a maximal torus of order $q^2-1=(2,q-1)\cdot k_1\cdot k_2$ implies the existence of an element of order  $rs$, i.~e., $r$ and $s$ are
adjacent in
$GK(G)$.
If $l=3$ (resp. $l=6$), then an element  $g$ of order $s$ is contained in a unique, up to conjugation, maximal torus of order $q^2+q+1=(3,q-1)k_3$
(resp. $q^2-q+1=(3,q+1)k_6$).
In this case $r,s$ are non-adjacent if and only if $r$ does not divide $\vert T\vert$, whence statement 1 of the lemma follows.

2. Since $\vert F_4(q)\vert=q^{24}(q^2-1)(q^6-1)(q^8-1)(q^{12}-1)$, the numbers $k,l$ are in the set $\{1,2,3,4,6,8,12\}$. If $l\le 3$, then the
existence  of maximal torus of order $(q^3-1)(q+1)=(2,q-1)\cdot(3,q-1)k_1\cdot k_2\cdot k_3$ implies that for every $k\le 3$ the primes $r,s$
are adjacent.  If $l=4$,
then an element    of order
$s$ of $G$ lies in a maximal torus $T$ of order equals to either
$(q-\epsilon)^2(q^2+1)$, or $(q^2-\epsilon)(q^2+1)$. In particular,
for every such maximal torus $T$ the inclusion $\pi(T)\subseteq R_1\cup R_2 \cup
R_4$ holds. Moreover there exists a maximal torus of order
$q^4-1=(2,q-1)^2\cdot k_1\cdot k_2\cdot k_4$.  So in
this case $r,s$ are non-adjacent if and only if $r$ does not divide $k_1\cdot k_2\cdot k_4$, i.~e., if and only if
$k=3$. If $l=6$, then each element  of order
$s$ of $G$ is in a maximal torus $T$ of order equals to either
$(q^3+1)(q-\epsilon)=(3,q+1)\cdot k_6\cdot (q+1)\cdot (q-\epsilon)$, or
$(q^2-q+1)^2=(3,q+1)^2\cdot k_6^2$.
In particular, for  $T$   the inclusion $\pi(T)\subseteq R_1\cup
R_2\cup R_6$ holds. Moreover there exists a maximal torus of order
$(q^3+1)(q-1)=(2,q-1)(3,q+1)k_1\cdot k_2\cdot k_6$. Thus $r,s$ are non-adjacent
if and only if  $k\in\{3,4\}$. If, finally, $l=8$ (resp. $l=12$),
then every element of order  $s$ of $G$ lies in a
maximal torus of order $(2,q-1)k_8$ (resp. $k_{12}$). Thus $r,s$ are
non-adjacent if and only if $k\not=8$ (resp.~${k\not=12}$).

3. Since $\vert E_6(q)\vert=\frac{1}{(3,q-1)}q^{36}(q^2-1)(q^5-1)(q^6-1)(q^8-1)(q^9-1)(q^{12}-1)$, the numbers $k,l$ are in the set
$\{1,2,3,4,5,6,8,9,12\}$. If $l\le 3$, then the existence in $G$ of a maximal
torus $T$ of order
$\frac{1}{(3,q-1)}(q^3-1)(q^2-1)(q-1)=(2,q-1)\cdot k_3\cdot k_2\cdot k_1\cdot(q-1)^2$
implies that $r,s$
are adjacent. If $l=4$, then each element of order $s$ of $G$ is in a maximal
torus of order equals either
$\frac{1}{(3,q-1)}(q^4-1)(q-\epsilon_1)(q-\epsilon_2)=\frac{1}{(3,q-1)}\cdot (2,q-1)^2\cdot k_1\cdot k_2\cdot k_4\cdot (q-\epsilon_1)\cdot(q-\epsilon_2)$, or
$\frac{1}{(3,q-1)}(q^3+1)(q^2+1)(q-1)=\frac{1}{(3,q-1)}\cdot (2,q-1)^2\cdot (3,q+1)\cdot k_6\cdot k_4\cdot k_2\cdot k_1$, or
$\frac{1}{(3,q-1)}(q^2+1)^2(q-1)^2=\frac{1}{(3,q-1)}\cdot (2,q-1)^2\cdot k_4^2\cdot (q-1)^2$.
Thus  $r,s$ are non-adjacent if and only if $k=3$. If $l=5$, then each element
of order $s$ of $G$ is in a maximal torus of order
$\frac{1}{(3,q-1)}(q^5-1)(q-\epsilon)=\frac{1}{(3,q-1)}(5,q-1)k_5(q-1)(q-\epsilon)$. Thus $r,s$ are non-adjacent if and only if $k\in\{3,4\}$.
If $l=6$, then every element of order $s$ of $G$ is in a maximal torus of order
equals either
$\frac{1}{(3,q-1)}(q^3+1)(q^2+q+1)(q-\epsilon)=(3,q+1)\cdot k_6\cdot k_3\cdot (q+1)\cdot (q-\epsilon)$, or
$\frac{1}{(3,q-1)}(q^3+1)(q^2+1)(q-1)=\frac{1}{(3,q-1)}\cdot(3,q+1)\cdot k_6\cdot (2,q-1)\cdot k_1\cdot k_2\cdot k_4$, or
$\frac{1}{(3,q-1)}(q^3+1)(q^2-1)(q-1)=\frac{1}{(3,q-1)}\cdot(3,q+1)\cdot k_6\cdot (2,q-1)^2\cdot k_1^2\cdot k_2^2$, or
$\frac{1}{(3,q-1)}(q^2+q+1)(q^2-q+1)^2=(3,q+1)^2\cdot k_6^2\cdot k_3$.
Thus $r,s$ are non-adjacent if and only if $k=5$. If $l=8$, then each element of
order $s$ of $G$ is in a
maximal torus of order
$\frac{1}{(3,q-1)}(q^4+1)(q^2-1)=\frac{1}{(3,q-1)}\cdot (2,q-1)^2\cdot k_8\cdot k_2\cdot k_1$. Hence  $r,s$ are non-adjacent if and only if either
$k\ge3$ and $k\not=8$, or
$r=3$ and $(q-1)_3=3$. If  $l=9$, then each element of order  $s$ of $G$ is in a
  maximal torus
of order $\frac{1}{(3,q-1)}(q^6+q^3+1)=k_9$. Hence $r,s$ are non-adjacent if and only if $k\not=9$. If, finally, $l=12$, then every element of order
$s$ of $G$ is in
maximal torus of order $\frac{1}{(3,q-1)}(q^4-q^2+1)(q^2+q+1)=k_{12}\cdot k_3$.
So $r,s$ are non-adjacent if and only if~${k\not=3,12}$.

4. Since $\vert {}^2E_6(q)\vert=\frac{1}{(3,q+1)}q^{36}(q^2-1)(q^5+1)(q^6-1)(q^8-1)(q^9+1)(q^{12}-1)$, the numbers $k,l$ are in the set
$\{1,2,3,4,6,8,10,12,18\}$. If $l\le 4$, the existence in $G$ of maximal tori of
orders
$\frac{1}{(3,q+1)}(q^3-1)(q^2+1)(q+1)=\frac{1}{(3,q+1)}\cdot (2,q-1)\cdot k_1\cdot k_2\cdot k_3\cdot k_4$,
$\frac{1}{(3,q+1)}(q^2+1)^2(q+1)^2=\frac{1}{(3,q+1)}\cdot(2,q-1)^2\cdot k_4^2\cdot (q+1)^2$, and
$\frac{1}{(3,q+1)}(q^3-1)(q^2-1)(q+1)=\frac{1}{(3,q+1)}\cdot (3,q-1)\cdot
(2,q-1)^2\cdot k_3\cdot k_1^2\cdot k_2^2$ implies that $r,s$ are adjacent.
If $l=6$, then each element of order $s$ of $G$ is contained in a maximal torus
of order equals either
$\frac{1}{(3,q+1)}(q^3+1)^2=(3,q+1)\cdot (q+1)^2\cdot k_6^2$, or $\frac{1}{(3,q+1)}(q^3+1)(q+1)(q-\epsilon_1)(q-\epsilon_2)=k_6(q+1)^2
(q-\epsilon_1)(q-\epsilon_2)$, or
$\frac{1}{(3,q+1)}(q^2-q+1)(q^3-\epsilon)(q-1)=(3,q+1)\cdot k_6\cdot(q^3-\epsilon)\cdot (q-1)$, or
$\frac{1}{(3,q+1)}(q^2-q+1)(q^2+q+1)^2=k_6\cdot (3,q-1)^2\cdot k_3^2$,
or $\frac{1}{(3,q+1)}(q^4-q^2+1)(q^2-q+1)=k_{12}\cdot k_6$.
Thus  $r,s$ are non-adjacent if and only if $k=4$. If $l=8$, then every element
of order $s$ of $G$ is in a
maximal torus of order $\frac{1}{(3,q+1)}(q^4+1)(q^2-1)=\frac{1}{(3,q+1)}\cdot (2,q-1)^2\cdot k_8\cdot k_2\cdot k_1$. So $r,s$ are
non-adjacent if and only if either
$k\ge3$ and $k\not=8$, or $r=3$ and $(q+1)_3=3$. If $l=10$, then each element of
order $s$ of $G$ is in a maximal torus of order equals
$\frac{1}{(3,q+1)}(q^5+1)(q-\epsilon)=\frac{1}{(3,q+1)}\cdot k_{10}\cdot (q+1)\cdot(q-\epsilon)$. Hence  $r,s$ are non-adjacent if and only if
$k\ge 3$, $k\not=10$.
If $l=12$, then every element of order $s$ of $G$ is contained in a   maximal
torus of order
$\frac{1}{(3,q+1)}(q^4-q^2+1)(q^2-q+1)=k_{12}\cdot k_6$. Therefore $r,s$ are
non-adjacent if and only if $k\not=6,12$. If, finally $l=18$, then each
element of order $s$ of $G$ is
contained in a
maximal torus of order $\frac{1}{(3,q+1)}(q^6-q^3+1)=k_{18}$. Hence $r,s$ are
non-adjacent if and only if~${k\not=18}$.

5. Since  $\vert E_7(q)\vert=\frac{1}{(2,q-1)}q^{63}(q^2-1)(q^6-1)(q^8-1)(q^{10}-1)(q^{12}-1)(q^{14}-1)(q^{18}-1)$, the numbers $k,l$ are in
$\{1,2,3,4,5,6,7,8,9,10,12,14,18\}$. There exist maximal tori of $G$ of orders
equal
$\frac{1}{(2,q-1)}(q^5-1)(q^2+q+1)=\frac{1}{(2,q-1)}\cdot(3,q-1)\cdot(5,q-1)\cdot(q-1)\cdot k_5\cdot k_3$,
 $\frac{1}{(2,q-1)}(q^4-1)(q^3-1)=(2,q-1)\cdot k_4\cdot k_1\cdot
k_2\cdot k_3\cdot (3,q-1)\cdot (q-1)$ and $\frac{1}{(2,q-1)}(q^5-1)(q^2-1)=k_5\cdot k_2\cdot k_1\cdot (5,q-1)\cdot (q-1)$, so for $l\le 5$
and
$(k,l)\not=(4,5)$
the numbers $r,s$ are adjacent. Since for $l=5$ every element of order $s$ of
$G$ is contained in a maximal torus of order either
$\frac{1}{(2,q-1)}(q^5-1)(q-1)(q-\epsilon)$, or $\frac{1}{(2,q-1)}(q^5-1)(q^2+q+1)=\frac{1}{(2,q-1)}\cdot(3,q-1)\cdot(5,q-1)\cdot(q-1)\cdot
k_5\cdot k_3$, we obtain that $r,s$ are non-adjacent if $(k,l)\not=(4,5)$. If
$l=6$, then the existence of maximal tori of $G$ of orders equal
$\frac{1}{(2,q-1)}(q^3+1)(q^4-1)=(2,q-1)\cdot(3,q+1)\cdot k_1\cdot k_2\cdot k_4 \cdot k_6\cdot (q+1)$ and
$\frac{1}{(2,q-1)}(q^6-1)(q-1)=(3,q^2-1) \cdot k_6\cdot k_3\cdot k_2\cdot k_1\cdot (q-1)$
implies that for $k\le 4$ and $k=6$ the numbers $r,s$ are adjacent.
Every element of order $s$ is in a maximal torus of order equals either
$\frac{1}{(2,q-1)}(q^3+1)(q^2+1)(q-\epsilon_1)(q-\epsilon_2)=(3,q+1)\cdot (q+1)\cdot k_6\cdot k_4\cdot (q-\epsilon_1)\cdot(q-\epsilon_2)$
with $(\epsilon_1,\epsilon_2)\not=(-1,-1)$, or
$\frac{1}{(2,q-1)}(q^3+1)(q-\epsilon_1)(q-\epsilon_2)(q-\epsilon_3)(q-\epsilon_4)=
\frac{1}{(2,q-1)}\cdot (3,q+1)\cdot k_6\cdot(q+1)(q-\epsilon_1)(q-\epsilon_2)(q-\epsilon_3)(q-\epsilon_4)$,
or $\frac{1}{(2,q-1)}(q^3+1)(q^3-\epsilon_1)(q-\epsilon_2)$, or
$\frac{1}{(2,q-1)}(q^2-q+1)^3(q+1)=\frac{1}{(2,q-1)}\cdot(3,q+1)^3\cdot k_6^3\cdot(q+1)$, or
$\frac{1}{(2,q-1)}(q^5+1)(q^2-q+1)=\frac{1}{(2,q-1)}\cdot(3,q+1)\cdot(5,q+1)\cdot k_6\cdot k_{10}\cdot(q+1)$,
or $\frac{1}{(2,q-1)}(q^3-1)(q^2-q+1)^2=\frac{1}{(2,q-1)}\cdot(q-1)\cdot (3,q-1)\cdot k_3\cdot (3,q+1)^2\cdot k_6^2$. Since for
$k=5$ the prime $r$ does not divide these numbers, we obtain that $r,s$ are non-adjacent if and only if $k=5$. If $l=7$, then each element
of order $s$ of $G$ is in a maximal torus of order
$\frac{1}{(2,q-1)}(q^7-1)=\frac{1}{(2,q-1)}\cdot (7,q-1)\cdot k_7\cdot(q-1)$.
Hence
$r,s$ are non-adjacent if and only if $k\not=1,7$. If $l=8$,  then every element
of order $s$ of $G$ is in a maximal torus of
order equals
$\frac{1}{(2,q-1)}(q^4+1)(q^2-\epsilon_1)(q-\epsilon_2)=k_8\cdot(q^2-\epsilon_1)(q-\epsilon_2)$. Hence $r,s$ are non-adjacent if and only if
$k\ge 3,k\not=4$. If $l=9$, then an element of order $s$ of $G$ is contained
in a   maximal torus of order
$\frac{1}{(2,q-1)}(q-1)(q^6+q^3+1)=\frac{1}{(2,q-1)}\cdot(q-1)\cdot(3,q-1)\cdot
k_9$.
Therefore $r,s$ are non-adjacent if and
only if $k\not=1,9$. If $l=10$, then an element of order $s$ of $G$ is contained
in a maximal torus of order equals either
$\frac{1}{(2,q-1)}(q^5+1)(q-1)(q-\epsilon)=(2,q-1)\cdot (5,q+1) \cdot k_{10}\cdot k_2\cdot k_1\cdot (q-\epsilon)$,  or
$\frac{1}{(2,q-1)}(q^5+1)(q^2-q+1)=\frac{1}{(2,q-1)}\cdot (5,q+1)\cdot (q+1)\cdot k_{10}\cdot(3,q+1)\cdot k_6$.
So $r,s$ are non-adjacent if and only if
$k\ge3$ and $k\not=6$. If $l=12$, then each element of order $s$ is contained in a maximal torus
of order equals $\frac{1}{(2,q-1)}(q^3-\epsilon)(q^4-q^2+1)=\frac{1}{(2,q-1)}\cdot (q^3-\epsilon)\cdot k_{12}$.
Hence  $r,s$ are non-adjacent if and only if $k\ge 4$ and $k\not=6,12$. If
$l=14$, then an element of order $s$ of $G$ is contained in
a maximal torus of order
$\frac{1}{(2,q-1)}(q^7+1)=\frac{1}{(2,q-1)}\cdot(7,q+1)\cdot k_{14}\cdot(q+1)$.
Therefore $r,s$ are non-adjacent if and only if $k\not=2,14$. If, finally, $l=18$,
then an element of order $s$ of $G$ is contained in a
maximal torus of order
$\frac{1}{(2,q-1)}(q+1)(q^6-q^3+1)=\frac{1}{(2,q-1)}\cdot(3,q+1)\cdot (q+1)\cdot k_{18}$. Therefore $r,s$ are non-adjacent if and only if~${k\not=2,18}$.

6. Since $\vert
E_8(q)\vert=q^{120}(q^2-1)(q^8-1)(q^{12}-1)(q^{14}-1)(q^{18}-1)(q^{20}-1)(q^{24}-1)(q^{30}-1)$, the
numbers $k,l$ are in the set $\{1,2,3,4,5,6,7,8,9,10,12,14,15,18,20,24,30\}$. Since $G$ contains
maximal tori of orders $(q^3-\epsilon_1)(q^4-1)(q-\epsilon_2)$,
$(q^5-1)(q^2+1)(q+1)=(5,q-1)\cdot k_5\cdot (2,q-1)^2\cdot k_4\cdot k_2\cdot k_1$, and
$(q^5-1)(q^3-1)=(3,q-1)\cdot (5,q-1)\cdot k_5\cdot k_3\cdot (q-1)^2$, for $l\le
6$ primes
$r,s$ are adjacent if $(k,l)\not=(5,6)$. If $k=5$, then every element of order
$r$ of $G$ is contained in a maximal torus of order equals
either $(q^5-1)(q^3-1)=(3,q-1)\cdot (5,q-1)\cdot k_5\cdot k_3\cdot (q-1)^2$, or $(q^5-1)(q^2+1)(q+1)=(5,q-1)\cdot k_5\cdot (2,q-1)^2\cdot
k_4\cdot k_2\cdot k_1$, or $(q^5-1)(q^2-1)(q-\epsilon)=(5,q-1)\cdot k_5\cdot (2,q-1)\cdot k_2\cdot k_1\cdot (q-1)\cdot (q-\epsilon)$, or
$(q^5-1)(q-1)^3$, or $(q^4+q^3+q^2+q+1)^2$, and all these orders are not divisible by $s$ for $l=6$. It follows that if $(k,l)=(5,6)$, then
$r,s$ are
non-adjacent. If $l=7$, then every element of order $s$ of $G$ is contained in a
maximal torus of
order $(q^7-1)(q-\epsilon)=(7,q-1)\cdot k_7\cdot (q-1)(q-\epsilon)$. So $r,s$ are non-adjacent if
and only if $k\ge3$ and $k\not=7$. If $l=8$, then an element  of order $s$ of
$G$ is contained in a
maximal torus of order equals either $(q^4+1)(q^4-\epsilon)=(2,q-1)\cdot k_8\cdot(q^4-\epsilon)$,
or $(q^4+1)(q^3-\epsilon_1)(q-\epsilon_2)=(2,q-1)\cdot k_8\cdot(q^3-\epsilon_1)(q-\epsilon_2)$ with $(\epsilon_1,\epsilon_2)\not=(-1,-1)$,
or
$(q^4+1)(q^2-1)^2=(2,q-1)\cdot k_8\cdot(q^2-1)^2$, or $(q^4+1)(q^2-\epsilon_1)
(q-\epsilon_2)^2=(2,q-1)\cdot k_8\cdot(q^2-\epsilon_1)\cdot (q-\epsilon_2)^2$. Hence $r,s$ are non-adjacent if and
only if $k=5,7$. If $l=9$, then an element of order $s$ of $G$ is contained in a
maximal torus of order
equals either $(q^6+q^3+1)(q-1)(q-\epsilon)=(3,q-1)\cdot k_9\cdot(q-1)\cdot(q-\epsilon)$, or
$(q^6+q^3+1)(q^2+q+1)=(3,q-1)^2\cdot k_9\cdot k_3$. Hence $r,s$ are non-adjacent if and only if
$k\ge4$ and $k\not=9$. If $l=10$ then every element of order $s$ of $G$ is
contained in a maximal torus of
order either $(q^5+1)(q^2-\epsilon_1)(q-\epsilon_2)=(5,q+1)\cdot k_{10}\cdot
(q+1)(q^2-\epsilon_1)(q-\epsilon_2)$ with $(\epsilon_1,\epsilon_2)\not=(-1,-1)$, or $(q^5+1)(q^3+1)=(5,q+1)\cdot k_{10}\cdot
(q+1)^2\cdot (3,q+1)\cdot k_6$, or $(q^5+1)(q^2-q+1)(q-1)=(5,q+1)\cdot k_{10}\cdot
(3,q+1)\cdot k_6\cdot(2,q-1)\cdot k_1\cdot k_2 $, or $(q^5+1)(q+1)^3=(5,q+1)\cdot k_{10}\cdot
(q+1)(q+1)^3$, or $(q^4-q^3+q^2-q+1)^2=((5,q+1)\cdot k_{10})^2$. Hence $r,s$ are
non-adjacent if and only if $k\ge3$ and $k\not=4,6,10$. If $l=12$,  then each
element of order $s$ of $G$ is
contained in a maximal torus of order equals either $(q^4-q^2+1)(q^2+1)(q^2-\epsilon)=(2,q-1)\cdot
k_{12}\cdot k_4\cdot(q^2-\epsilon)$, or $(q^4-q^2+1)(q^2+1)(q-\epsilon)^2=(2,q-1)\cdot k_{12}\cdot
k_4\cdot(q-\epsilon)^2$, or $(q^4-q^2+1)(q^3-\epsilon_1) (q-\epsilon_2)=k_{12}\cdot
(q^3-\epsilon_1)\cdot (q-\epsilon_2)$, or $(q^4-q^2+1)(q^2+q+1)^2=(3,q-1)^2\cdot k_{12}\cdot
k_3^2$, or $(q^4-q^2+1)(q^2-q+1)^2=(3,q+1)^2\cdot k_{12}\cdot k_6^2$. Hence $r,s$ are non-adjacent
if and only if  $k\ge5$ and $k\not=6,12$. If $l=14$, then an element of order
$s$ of $G$ is contained in a
maximal torus of order $(q^7+1)(q-\epsilon)=(7,q+1)\cdot k_{14}\cdot(q+1)\cdot(q-\epsilon)$.
Therefore $r,s$ are non-adjacent if and only if $k\ge3$ and $k\not=14$. If  $l=15,24,30$, then each
element of order $s$ of $G$ is contained in a   maximal
torus of order $k_{l}$. So
$r,s$ are non-adjacent if and only if $k\not=l$. If $l=18$, then an element of
order $s$ of $G$ is
contained in a maximal torus of order equals either $(q^6-q^3+1)(q+1)(q-\epsilon)=(3,q+1)\cdot
k_{18}\cdot(q+1)\cdot(q-\epsilon)$, or $(q^6-q^3+1)(q^2-q+1)=(3,q+1)^2\cdot k_{18}\cdot k_6$. Hence
$r,s$ are non-adjacent if and only if $k\ge 3$ and $k\not=6,18$. If $l=20$, then every element of
order $s$ of $G$ is contained in a  maximal torus of
order
$q^8-q^6+q^4-q^2+1=(5,q^2+1)\cdot k_{20}$. So $r,s$ are non-adjacent if and only if  $r\cdot
k\not=20$ (i.~e., $r\not=5$ or $k\not=4$) and $k\not=20$.

7. Since  $\vert {}^3D_4(q)\vert=q^{12}(q^2-1)(q^6-1)(q^8+q^4+1)$, the numbers  $k,l$ are in the set $\{1,2,3,6,12\}$. Since
$G$ contains maximal tori of orders $(q^3-\epsilon_1)(q-\epsilon_2)$, then for $l\le 3$ primes $r,s$ are adjacent. If $l=6$, then each element of
order $s$ of $G$ is in a maximal torus of order
$(q^3+1)(q-\epsilon)=(3,q+1)\cdot k_6\cdot (q+1)\cdot (q-\epsilon)$.
Hence $r,s$ are non-adjacent if and only if $k=3$.
 If $l=12$, then and
element of order $s$ of $G$ is contained in a   maximal
torus of order $q^4-q^2+1=k_{12}$ and $r,s$ are non-adjacent if and only
if~${k\not=12}$.
\end{proof}

Now we consider simple Suzuki and Ree groups.

\begin{lem}\label{SuzReeDivisors} Let $n$ be a natural number.

\noindent {\em 1.} Let
 $m_1(B,n)=2^{2n+1}-1$,

$m_2(B,n)=2^{2n+1}-2^{n+1}+1$,

$m_3(B,n)=2^{2n+1}+2^{n+1}+1$.

Then $(m_i(B,n),m_j(B,n))=1$ if~$i\not=j$.
%\item[{\em 2.}]

\noindent {\em 2.} Let $m_1(G,n)=3^{2n+1}-1$,

$m_2(G,n)=3^{2n+1}+1$,

$m_3(G,n)=3^{2n+1}-3^{n+1}+1$,

$m_4(G,n)=3^{2n+1}+3^{n+1}+1$.

Then $(m_1(G,n), m_2(G,n))=2$ and $(m_i(G,n),m_j(G,n))=1$ otherwise.
%\item[{\em 3.}]

\noindent {\em 3.} Let $m_1(F,n)=2^{2n+1}-1$,

$m_2(F,n)=2^{2n+1}+1$,

$m_3(F,n)=2^{4n+2}+1$,

$m_4(F,n)=2^{4n+2}-2^{2n+1}+1$,

$m_5(F,n)=2^{4n+2}-2^{3n+2}+2^{2n+1}-2^{n+1}+1$,

$m_6(F,n)=2^{4n+2}+2^{3n+2}+2^{2n+1}+2^{n+1}+1$.

Then $(m_2(F,n),m_4(F,n))=3$ and $(m_i(F,n),m_j(F,n))=1$ otherwise.
\end{lem}

\begin{proof} Items (1) and (2) are repeated items (1) and (2) of \cite[Lemma~1.5]{VasVd}. Item (3) is
corrected with respect to Lemma~\ref{toriofexcptgrps}.
\end{proof}

If $G$ is a Suzuki or a Ree group over a field of order $q$, then denote by $S_i(G)$ the set $\pi(m_i(B,n))$ for
$G={}^2B_2(2^{2n+1})$, the set $\pi(m_i(G,n))\setminus\{2\}$ for $G={}^2G_2(3^{2n+1})$, and the set
$\pi(m_i(F,n))\setminus\{3\}$ for $G={}^2F_4(2^{2n+1})$. If $G$ is fixed, then we put $S_i=S_i(G)$,
and denote by $s_i$ any prime from~$S_i$.

\begin{prop}\label{adjsuzree}
Let $G$ be a finite simple Suzuki or Ree group over a field of characte\-ristic~$p$, let $r,s$  be
odd primes with $r,s\in\pi(G)\setminus\{p\}$. Then $r,s$ are non-adjacent if and only if one of the
following holds:
\begin{itemize}
\item[{\em 1.}] $G={^2B_2(2^{2n+1})}$, $r\in S_k(G)$, $s\in S_l(G)$ and~$k\not=l$.
\item[{\em 2.}] $G={^2G_2(3^{2n+1})}$, $r\in S_k(G)$, $s\in S_l(G)$ and~$k\not=l$.
\item[{\em 3.}] $G={^2F_4(2^{2n+1})}$, either $r\in S_k(G)$, $s\in S_l(G)$ and~$k\not=l$,
$\{k,l\}\neq\{1,2\},\{1,3\}$; or $r=3$ and $s\in S_l(G)$, where $l\in\{3,5,6\}$.
\end{itemize}
\end{prop}

\begin{proof}
Follows from \cite[Lemma~1.3]{VasVd}, Lemma \ref{toriofexcptgrps}, and Lemma \ref{SuzReeDivisors}.
\end{proof}

\section{Cocliques for groups of Lie type}
Let $G$ be a finite simple group of Lie type with the base field of order $q$ and
characteristic~$p$. Every $r\in\pi(G)\setminus\{p\}$ is known to be a primitive prime divisor of
$q^i-1$, where $i$ is bounded by some function depending on the Lie rank of~$G$. Given a finite
simple group of Lie type $G$, define a set $I(G)$ as follows. If $G$ is neither a Suzuki, nor a Ree
group, then $i\in I(G)$ if and only if $\pi(G)\cap R_i(q)\not=\varnothing$. If $G$ is either a
Suzuki or a Ree group, then $i\in I(G)$ if and only if $\pi(G)\cap S_i(G)\not=\varnothing$. Notice
that if $\pi(G)\cap R_i(q)\not=\varnothing$ (resp.  $\pi(G)\cap S_i(G)\not=\varnothing$), then
$R_i(q)\subseteq\pi(G)$ (resp. $S_i(G)\subseteq\pi(G)$). Thus, the following partition of $\pi(G)$
arises:
$$\pi(G)=\{p\}\sqcup\bigsqcup_{i\in I(G)}R_i,$$
or $$\pi(G)=\{2\}\sqcup\bigsqcup_{i\in I(G)}S_i$$ in case of Suzuki groups, or
$$\pi(G)=\{2\}\sqcup\{3\}\cup\bigsqcup_{i\in I(G)}S_i$$ in case of Ree groups.

As followed from an adjacency criterion, two distinct primes from the same class of the partition
are always adjacent. Moreover, in most cases an answer to the question: whether two primes from
distinct classes $R_i$ and $R_j$ (or $S_i$ and $S_j$) of the partition are adjacent, depends only
on the choice of the indices $i$ and $j$. We formalize this inference by the following
definitions.

\begin{df}\label{MG} Suppose  $G$ is a finite simple group of Lie type with the base field of order $q$ and characteristic $p$, and $G$
is not
isomorphic to ${}^2B_2(2^{2m+1})$, ${}^2G_2(3^{2m+1})$, ${}^2F_4(2^{2m+1})$, and
$A^\varepsilon_2(q)$. Define the set $M(G)$ to be a subset of $I(G)$ such that
$i\in M(G)$ if and only if the
intersection of $R_i$ and every coclique of maximal size of $GK(G)$ is nonempty.
\end{df}

\begin{df}\label{MGSR} If $G={}^2B_2(2^{2m+1})$ or ${}^2G_2(3^{2m+1})$, $m\geqslant1$, then put $M(G)=I(G)$. If
$G={}^2F_4(2^{2m+1})$, $m\geqslant2$, then put $M(G)=\{2,3,4,5,6\}$. If $G={}^2F_4(8)$, then put
$M(G)=\{5,6\}$.
\end{df}

\begin{df}\label{theta}
Suppose $G$ is a finite simple group of Lie type with the base field of order $q$ and characteristic $p$, and $G$ is not
isomorphic to  ${}^2B_2(2^{2m+1})$, ${}^2G_2(3^{2m+1})$, ${}^2F_4(2^{2m+1})$, and
$A^\varepsilon_2(q)$. A set $\Theta(G)$ consists of all subsets $\theta(G)$ of $\pi(G)$ satisfying
the following conditions:

(a) $p$ lies in $\theta(G)$ if and only if $p$ lies in every coclique of maximal size of $GK(G)$;

(b) for every $i\in M(G)$ exactly one prime from $R_i$ lies in $\theta(G)$.
\end{df}

\begin{df}\label{thetaS}
Let $G={}^2B_2(2^{2m+1})$. A set $\Theta(G)$ consists of all subsets $\theta(G)$ of $\pi(G)$
satisfying the following conditions:

(a) $p=2$ lies in $\theta(G);$

(b) for every $i\in M(G)$ exactly one prime from $S_i$ lies in $\theta(G)$.
\end{df}

\begin{df}\label{thetaR3}
Let $G={}^2G_2(3^{2m+1})$. A set $\Theta(G)$ consists of all subsets $\theta(G)$ of $\pi(G)$
satisfying the following conditions:

(a) $p=3$ lies in $\theta(G);$

(b) for every $i\in M(G)$ exactly one prime from $S_i$ lies in $\theta(G).$
\end{df}

\begin{df}\label{thetaR2}
Let $G={}^2F_4(2^{2m+1})$, $m\geqslant1$. A set $\Theta(G)$ consists of all subsets $\theta(G)$ of
$\pi(G)$ satisfying the following condition:

(a) for every $i\in M(G)$ exactly one prime from $S_i$ lies in $\theta(G).$

\end{df}

\begin{df}\label{theta2A} Let $G=A^{\varepsilon}_2(q)$, and $(q,\varepsilon)\neq(2,-)$. If $q+\varepsilon1\neq2^k$, then put
$M(G)=\{\nu_{\varepsilon}(2),\nu_{\varepsilon}(3)\}$, and if $q+\varepsilon1=2^k$, then
$M(G)=\{\nu_{\varepsilon}(3)\}$. A set $\Theta(G)$ consists of all subsets $\theta(G)$ of $\pi(G)$
satisfying the following conditions.

(1) $p$ lies in $\theta(G)$ if and only if $q+\varepsilon1\neq2^k;$

(2) if $(q-\varepsilon1)_3=3$, then $3\in\theta(G)$.

(3) for every $i\in M(G)$ exactly one prime from $R_{\nu_{\varepsilon}(i)}$ lies in $\theta(G)$,
excepting one case: if $2\in R_{\nu_{\varepsilon}(2)}$, then $2$ does not lie in $\theta(G)$.
\end{df}

\textsl{Remark.} A function $\nu_{\varepsilon}$ is defined in~\eqref{nuepsilon(n)}.

\begin{df}\label{Thetaprime} Let $G$ be a finite simple group of Lie type. The subset $\theta'(G)$
of $\pi(G)$ is an element of $\Theta'(G)$, if for every $\theta(G)\in\Theta(G)$ the union
$\rho(G)=\theta(G)\cup\theta'(G)$ is a coclique of maximal size in $GK(G)$.
\end{df}

Now we describe cocliques of maximal size for groups of Lie type. First we consider classical
groups postponing groups $A_1(q)$, $A_2^{\varepsilon}(q)$ to the end of the section.

\begin{prop}\label{cocliqueAE4} If $G$ is one of finite simple groups $A_{n-1}(q)$, ${}^2A_{n-1}(q)$
with the base field of characteristic~$p$ and order $q$, and $n\geqslant4$, then $t(G)$, and the sets $\Theta(G)$,
$\Theta'(G)$ are listed in Table~{\em\ref{LinearUnitaryTable}}.
\end{prop}

\begin{proof}

It is obvious that the function $\nu_{\varepsilon}$ defined in~\eqref{nuepsilon(n)} is a bijection on $\mathbb{N}$, so
$\nu_{\varepsilon}^{-1}$
is well defined. Moreover, since $\nu_{\varepsilon}^2$ is the identity map, we have
$\nu_{\varepsilon}^{-1}=\nu_{\varepsilon}$.

Using Lemma \ref{Zsigmondy Theorem} and an information on the orders of groups
$A^{\varepsilon}_{n-1}(q)$, we obtain that a number $i$ lies in $I(G)$, if the following conditions
holds:

(a) $\nu_{\varepsilon}(i)\leqslant n$;

(b) $i\neq1$ for $q=2,3$, and $i\neq6$ for $q=2$.

By \cite[Propositions~2.1, 2.2, 4.1, and~4.2]{VasVd} two distinct primes from $R_i$ are adjacent for every~${i\in I(G)}$.

Denote by $N(G)$ the set $\{i\in I(G)\mid n/2<\nu_{\varepsilon}(i)\leqslant n\}$ and by $\chi$ any
set of type $\{r_i\mid i\in N(G)\}$ such that $|\chi\cap R_i|=1$ for all $i\in N(G)$. Note that
$1,2$ can not lie in $N(G)$, because $n\geqslant4$. In particular, $2$ does not lie in any $\chi$.
Let $i\neq j$ and $n/2<\nu_{\varepsilon}(i),\nu_{\varepsilon}(j)\leqslant n$. Then
$\nu_{\varepsilon}(i)+\nu_{\varepsilon}(j)>n$ and $\nu_{\varepsilon}(i)$ does not divide
$\nu_{\varepsilon}(j)$. By \cite[Propositions~2.1,~2.2]{VasVd}, primes $r_i$ and $r_j$ are not
adjacent. Thus, every $\chi$ forms a coclique of $GK(G)$.

Denote by $\xi$ the set
$$\{p\}\cup\bigcup_{i\in I(G)\setminus N(G)}R_i.$$ By
\cite[Propositions 2.1, 2.2, 3.1, 4.1, 4.2]{VasVd} every two distinct primes from $\xi$ are
adjacent in $GK(G)$. Thus, every coclique of $GK(G)$ contains at most one prime from~$\xi$.

\textbf{Case 1.} Let $n\geqslant7$.

If $q=2$ and $G=A_{n-1}(q)$ we assume that $n\geqslant13$ first, in order to avoid the exceptions
arising because of $R_6=\varnothing$ for $q=2$.

The conditions on $n$ implies that $|N(G)|\geqslant4$. By \cite[Proposition~3.1]{VasVd}, we have
that $t(p,G)\leqslant3$, so $p$ can not lie in any coclique of maximal size. By
\cite[Propositions~4.1, 4.2]{VasVd}, the same assertion is true for any primitive prime divisor
$r_i$, where $\nu_{\varepsilon}(i)=1$. Thus, solving the problem, does a prime $r$ lie in a
coclique of maximal size of $GK(G)$, we may assume that $r$ is neither a characteristic, nor a
divisor of $q-\varepsilon1$. Hence \cite[Propositions~2.1 and~2.2]{VasVd} will be the main
technical tools.

Suppose that $n=2t+1$ is odd. If $\nu_{\varepsilon}(i)\leqslant n/2$, then there exist at least two
distinct numbers $j,k$ from $N(G)$ such that $r_i$ is adjacent to $r_j$ and $r_k$. Indeed, if
$\nu_{\varepsilon}(i)<t$, then we take $j$ and $k$ such that $\nu_{\varepsilon}(j)=t+1$ and
$\nu_{\varepsilon}(k)=t+2$, while if $\nu_{\varepsilon}(i)=t$, then we take $j$ and $k$ such that
$\nu_{\varepsilon}(j)=t+1$ and $\nu_{\varepsilon}(k)=2t$. Thus, $M(G)=N(G)$, every
$\theta(G)\in\Theta(G)$ is of type $\{r_i\mid n/2<\nu_{\varepsilon}(i)\leqslant n\}$,
$\Theta'(G)=\varnothing$, and $t(G)=t+1=[(n+1)/2]$.

Suppose that $n=2t$ is even. If $\nu_{\varepsilon}(i)< n/2$, then there exist at least two distinct
numbers $j,k$ from $N(G)$ such that $r_i$ is adjacent to $r_j$ and $r_k$. It is sufficient, to take
$j$ and $k$ such that $\nu_{\varepsilon}(j)=t+1$, and $\nu_{\varepsilon}(k)=t+2$ if
$\nu_{\varepsilon}(i)<t-1$, or $\nu_{\varepsilon}(k)=2t-2$ if $\nu_{\varepsilon}(i)=t-1$. On the
other hand, if $\nu_{\varepsilon}(i)=t=n/2$, then $r_i$ is adjacent to $r_j$, where
$\nu_{\varepsilon}(j)=2t=n$, and is non-adjacent to every $r_k$, where $k\in N(G)$ and $k\neq j$.
Thus, $M(G)=N(G)\setminus\{\nu_{\varepsilon}(n)\}$, every $\theta(G)\in\Theta(G)$ is of type
$\{r_i\mid n/2<\nu_{\varepsilon}(i)<n\}$, and $\Theta'(G)$ consists of one-element sets of type
$\{r_{\nu_{\varepsilon}(n/2)}\}$ or $\{r_{\nu_{\varepsilon}(n)}\}$. Hence, $t(G)=t=[(n+1)/2]$.

It remains to consider the following  cases: $q=2$, $G=A_{n-1}(q)$, and
$7\leqslant n\leqslant12$. All
results (see Table~\ref{LinearUnitaryTable}) are obtained by arguments similar to that in general
case with respect to the fact: $R_6=\varnothing$, and can be easily verified by using
\cite[Propositions 2.1, 2.2, 3.1, 4.1, 4.2]{VasVd}. The most interesting case arises, when $n=8$.
In that case $\Theta(G)$ consists of one-element sets $\theta(G)$ of type $\{r_7\}$, while
$\Theta'(G)$ consists of two-elements sets $\theta'(G)$ of types $\{p,r_8\}$, $\{r_4,r_5\}$, $\{r_3,r_8\}$, or
$\{r_5,r_8\}$.

\textbf{Case 2.} Let $n=6$.

First, we assume that $q\neq2$. Then
$N(G)=\{\nu_{\varepsilon}(4),\nu_{\varepsilon}(5),\nu_{\varepsilon}(6)\}$, and $|N(G)|=3$.
Therefore, a set of type
$\{r_{\nu_{\varepsilon}(4)},r_{\nu_{\varepsilon}(5)},r_{\nu_{\varepsilon}(6)}\}$ forms a coclique
in $GK(G)$, and $t(G)\geqslant3$. Arguing as in previous case, we obtain that any prime
$r_{\nu_{\varepsilon}(3)}$ is adjacent to $r_{\nu_{\varepsilon}(6)}$, and a set of type
$\{r_{\nu_{\varepsilon}(3)},r_{\nu_{\varepsilon}(4)},r_{\nu_{\varepsilon}(5)}\}$ is a coclique. By
\cite[Proposition~3.1]{VasVd}, we have that a set of type
$\{p,r_{\nu_{\varepsilon}(5)},r_{\nu_{\varepsilon}(6)}\}$ is a coclique, and $p$ is adjacent to any
prime $r_{\nu_{\varepsilon}(4)}$. If $\nu_{\varepsilon}(i)=2$, then $r_i$ is adjacent to $p$, and
is non-adjacent to $r_j$ if and only if $\nu_{\varepsilon}(j)=5$, so $t(r_i,G)=2$. Let $r$ be a
divisor of $q-\varepsilon1$. If $r\neq3$ or $(q-\varepsilon1)_3\neq3$, then \cite[Propositions~4.1,
4.2]{VasVd} implies that $t(r,G)=2$, while if $r=3$ and $(q-\varepsilon1)_3=3$, we have that
$t(3,G)=3$ and a set of type $\{3,r_{\nu_{\varepsilon}(5)},r_{\nu_{\varepsilon}(6)}\}$ is a
coclique in $GK(G)$. Thus, if $q\neq2$, then $M(G)=\{\nu_{\varepsilon}(5)\}$, and every
$\theta(G)\in\Theta(G)$ is of type $\{r_{\nu_{\varepsilon}(5)}\}$. Every $\theta'(G)\in\Theta'(G)$
is a two-element set of type $\{p,r_{\nu_{\varepsilon}(6)}\}$,
$\{r_{\nu_{\varepsilon}(3)},r_{\nu_{\varepsilon}(4)}\}$,
$\{r_{\nu_{\varepsilon}(4)},r_{\nu_{\varepsilon}(6)}\}$, and if $(q-\varepsilon1)_3=3$ is also of
type $\{3,r_{\nu_{\varepsilon}(6)}\}$.

Let $G=A_5(2)$. Since $R_6=R_1=\varnothing$, we have that every $\theta(G)\in\Theta(G)$ is of type
$\{r_{3},r_{4},r_{5}\}$, and $\Theta'(G)=\varnothing$.

Let $G={}^2A_5(2)$. Since $R_6=R_{\nu(3)}=\varnothing$ and $(q+1)_3=3$, we have that every
$\theta(G)\in\Theta(G)$ is of type $\{r_{3},r_{10}\}$, and every $\theta'(G)\in\Theta'(G)$ is a one
element set of type $\{p\}$, $\{r_4\}$, or $\{3\}$.

In all cases $t(G)=3$.

\textbf{Case 3.} Let $n=5$.

We have $N(G)=\{\nu_{\varepsilon}(4),\nu_{\varepsilon}(5)\}$, and $|N(G)|=2$, so $t(G)\leqslant3$.
Assume now that $G\neq{}^2A_{4}(2)$. Then $R_{\nu_{\varepsilon}(3)}$ is always nonempty, and a set
of type $\{r_{\nu_{\varepsilon}(3)},r_{\nu_{\varepsilon}(4)},r_{\nu_{\varepsilon}(5)}\}$ is a
coclique in $GK(G)$. By \cite[Proposition~3.1]{VasVd}, we have that a set of type
$\{p,r_{\nu_{\varepsilon}(4)},r_{\nu_{\varepsilon}(5)}\}$ is also a coclique. A prime
$r_{\nu_{\varepsilon}(2)}$ is adjacent to $p$, and is non-adjacent to $r_j$ if and only if
$\nu_{\varepsilon}(j)=5$. Let $r$ be a divisor of $q-\varepsilon1$. If $r\neq5$ or
$(q-\varepsilon1)_5\neq5$, then \cite[Propositions~4.1, 4.2]{VasVd} implies that $t(r,G)=2$, while
if $r=5$ and $(q-\varepsilon1)_5=5$, we have that $t(5,G)=3$ and a set of type
$\{5,r_{\nu_{\varepsilon}(4)},r_{\nu_{\varepsilon}(5)}\}$ is a coclique in $GK(G)$. Thus, if
$G\neq{}^2A_{4}(2)$, then $M(G)=N(G)$, every $\theta(G)\in\Theta(G)$ is of type
$\{r_{\nu_{\varepsilon}(4)},r_{\nu_{\varepsilon}(5)}\}$. Every $\theta'(G)\in\Theta'(G)$ is a
one-element set of type $\{p\}$ or $\{r_{\nu_{\varepsilon}(3)}\}$, and if $(q-\varepsilon1)_5=5$ is
also of type $\{5\}$ .

Let $G={}^2A_{4}(2)$. Since $R_6=R_{\nu(3)}=\varnothing$ and $(q+1)_5=1$, we have that every
$\theta(G)\in\Theta(G)$ is of type $\{p,r_{4},r_{10}\}$, and $\Theta'(G)=\varnothing$.

In all cases $t(G)=3$.

\textbf{Case 4.} Let $n=4$.

First, we assume that $G\neq{}^2A_{3}(2)$. Then
$N(G)=\{\nu_{\varepsilon}(3),\nu_{\varepsilon}(4)\}$, and $|N(G)|=2$, so $t(G)\leqslant3$. By
\cite[Proposition~3.1]{VasVd}, we have that a set of type
$\{p,r_{\nu_{\varepsilon}(3)},r_{\nu_{\varepsilon}(4)}\}$ is a coclique in $GK(G)$. A prime
$r_{\nu_{\varepsilon}(2)}$ is adjacent to $p$ and any prime $r_{\nu_{\varepsilon}(4)}$. If $r$ is
an odd prime divisor of $q-\varepsilon1$, then \cite[Propositions~4.1, 4.2]{VasVd} implies that
$t(r,G)=2$. The same assertion is true for $r=2$ if and only if $(q-\varepsilon1)_2\neq4$, while if
$(q-\varepsilon1)_2=4$ then $\{2,r_{\nu_{\varepsilon}(3)},r_{\nu_{\varepsilon}(4)}\}$ is a
coclique. Therefore, if $(q-\varepsilon1)_2\neq4$, then every $\theta(G)\in\Theta(G)$ is of type
$\{p,r_{\nu_{\varepsilon}(3)},r_{\nu_{\varepsilon}(4)}\}$, and $\Theta'(G)=\varnothing$. But if
$(q-\varepsilon1)_2=4$, then $M(G)=N(G)$, every $\theta(G)\in\Theta(G)$ is of type
$\{r_{\nu_{\varepsilon}(3)},r_{\nu_{\varepsilon}(4)}\}$, and $\Theta'(G)=\{\{2\},\{p\}\}$. Anyway,
$t(G)=3$.

Let $G={}^2A_{3}(2)$. Since $R_6=R_{\nu(3)}=\varnothing$ and $(q+1)_4=1$, we obtain that
$\Theta(G)=\{\{r_4\}\}$, and $\Theta'(G)=\{\{p\},\{r_2\}\}$. In this case, $t(G)=2$.
\end{proof}

\begin{prop}\label{cocliqueBC} If $G$ is one of finite simple groups $B_n(q)$, $C_n(q)$, $D_n(q)$ or ${}^2D_n(q)$
with the base field of characteristic~$p$ and order $q$, then $t(G)$, and the sets $\Theta(G)$, $\Theta'(G)$ are listed
in Table~{\em\ref{ClassicTable}}.
\end{prop}

\begin{proof} Using Lemma \ref{Zsigmondy Theorem} and an information on the
orders of groups under consideration, we obtain that a number $i$ lies in $I(G)$, if the following
conditions holds:

(a) $\eta(i)\leqslant n$;

(b) $i\neq1$ for $q=2,3$, and $i\neq6$ for $q=2$;

(c) $i\neq2n$ for $G=D_n(q)$;

(d) $i\neq n$ for $G={}^2D_n(q)$ and $n$ odd.

By \cite[Proposition~4.3 and~4.4]{VasVd} and Propositions \ref{adjbn} and \ref{adjdn} it follows that for every $i\in I(G)$ two distinct
primes from $R_i$ are adjacent.

Denote by $N(G)$ the set $\{i\in I(G)\mid n/2<\eta(i)\leqslant n\}$ and by $\chi$ any set of type
$\{r_i\mid i\in N(G)\}$ such that $|\chi\cap R_i|=1$ for all $i\in N(G)$. Let $i\neq j$ and
$n/2<\eta(i),\eta(j)\leqslant n$. We have $\eta(i)+\eta(j)>n$. Suppose that $i/j$ is an odd
natural number. Then $i$ and $j$ is of the same parity, so $\eta(i)/\eta(j)$ is also an odd
natural number. Since $i\neq j$, we have $\eta(j)>2\eta(i)>n$, contrary to the choice of~$j$. Thus
$i/j$ is not an odd number. By Propositions~\ref{adjbn} and~\ref{adjdn}, primes $r_i$ and $r_j$
are not adjacent. Thus, every $\chi$ forms a coclique of $GK(G)$.

Denote by $\xi$ the set
$$\{p\}\cup\bigcup_{i\in I(G)\setminus N(G)}R_i.$$ By Propositions~\ref{adjbn},~\ref{adjdn} and
\cite[Propositions 3.1, 4.3, 4.4]{VasVd} every two distinct primes from $\xi$ are adjacent in
$GK(G)$. Thus, every coclique of $GK(G)$ contains at most one prime from~$\xi$.

Now we determine cocliques of maximal size considering the groups of different types separately.
However, by \cite[Theorem~7.5]{VasVd}, we have $GK(B_n(q))=GK(C_n(q))$, and so analysis for groups
of types $B_n$ and $C_n$ is mutual.

\textbf{Case 1.} Let $G$ be one of the simple groups $B_n(q)$ or $C_n(q)$.

Suppose that $n=2$. If $q=2$, then the group $G$ is not simple, so we can assume that
$q\geqslant3$. If $q=3$, then $I(G)=\{2,4\}$, and if $q>3$, then $I(G)=\{1,2,4\}$. In both cases,
$N(G)=\{4\}$. Since $r_4$ is non-adjacent to every $r\in\xi$, we have $M(G)=N(G)=\{4\}$, every
$\theta(G)\in\Theta(G)$ is a one-element set containing exactly one element $r_4$ from $R_4$.
Every $\theta'(G)\in\Theta'(G)$ is a one-element set containing exactly one element from~$\xi$.
Thus, $t(G)=2$.

Suppose that $n=3$. If $q\neq2$ then $N(G)=\{3,6\}$, and if $q=2$ then $N(G)=\{3\}$. The set
$\{1,2,3,4,6\}$ includes $I(G)$, and so $\xi=\{p\}\cup R_1\cup R_2\cup R_4$, where $\{p\}$, $R_2$ and
$R_4$ are always nonempty. The prime $p$ and any prime $r_4$ are adjacent one to another, and are
non-adjacent to every $r_i$ with $i\in N(G)$. On the other hand, for $i\in\{1,2\}$ and
$j\in\{3,6\}$, primes $r_i$ and $r_j$ are adjacent. Therefore, $M(G)=N(G)$, $\theta(G)$ is
of type $\{r_3\}$ for $q=2$, and is of type $\{r_3,r_6\}$ otherwise. The set $\Theta'(G)$ consists
of one-element sets of type $\{p\}$, $\{r_2\}$, and $\{r_4\}$, if $q=2$, and sets of type $\{p\}$,
and $\{r_4\}$ otherwise. Thus, $t(G)=2$ for $q=2$, and $t(G)=3$ otherwise.

Let $n\geqslant 4$. Now we consider four different cases subject to residue of $n$ modulo $4$. We
write $n=4t+k$, where $k=0,1,2,3$, and $t\geqslant1$. If $q=2$ we assume that $t>1$ to avoid
exceptional cases that arise because of $R_6=\varnothing$ for $q=2$.

Suppose that $n=4t$. Then
$$N(G)=\{2t+1,2t+3,\ldots,4t-1,4t+2,4t+4,\ldots,8t\},$$
and so $|N(G)|=3t$. By adjacency criterion, $r_{4t}$ is non-adjacent to every $r_i$, where $i\in
N(G)$. Therefore, $t(G)\geqslant3t+1\ge4$. By \cite[Propositions 3.1, 4.3]{VasVd}, we have
$t(2,G)\leqslant t(p,G)<4$, so $p$ and $2$ cannot lie in any coclique of maximal size. Furthermore,
if $\eta(i)<n/2=2t$, then any odd $r_i$ is adjacent to $r_{4t}$, $r_{2t+1}$, $r_{4t+2}$. Therefore,
$M(G)=N(G)\cup\{n\}$, every $\theta(G)\in\Theta(G)$ is of type $\{r_i\mid
n/2\leqslant\eta(i)\leqslant n\}$, $\Theta'(G)=\varnothing$, so $t(G)=3t+1=[(3n+5)/4]$.

Suppose that $n=4t+1$. Then
$$N(G)=\{2t+1,2t+3,\ldots,4t+1,4t+2,4t+4,\ldots,8t+2\},$$
so $|N(G)|=3t+2$ and $t(G)\geqslant5$. By \cite[Propositions 3.1, 4.3]{VasVd}, we have
$t(2,G)\leqslant t(p,G)<4$. Therefore, $p$ and $2$ cannot lie in any coclique of maximal size. If
$\eta(i)<n/2$, then any odd $r_i$ is adjacent to $r_{2t+1}$, $r_{4t+2}$, so cannot lie in any
coclique of maximal size. Thus, $M(G)=N(G)$, every $\theta(G)\in\Theta(G)$ is of type $\{r_i\mid
n/2<\eta(i)\leqslant n\}=\{r_i\mid n/2\leqslant\eta(i)\leqslant n\}$, $\Theta'(G)=\varnothing$, and
$t(G)=3t+2=[(3n+5)/4]$.

Suppose that $n=4t+2$. Then
$$N(G)=\{2t+3,2t+5,\ldots,4t+1,4t+4,4t+6,\ldots,8t+4\},$$
so $|N(G)|=3t+1$ and $t(G)\geqslant4$. Since $t(2,G)\leqslant t(p,G)<4$, primes $p$ and $2$ cannot
lie in any coclique of maximal size. Any primes $r_{2t+1}$ and $r_{4t+2}$ are adjacent one to
another and are non-adjacent to every $r_i$ with $i\in N(G)$. If $\eta(i)<n/2$, then $r_i$ is
adjacent to $r_{4t+4}$, $r_{4t+2}$, and $r_{2t+1}$. Therefore, $N(G)=M(G)$, every
$\theta(G)\in\Theta(G)$ is of type $\{r_i\mid n/2<\eta(i)\leqslant n\}$, and $\Theta'(G)$ consists
of one-element sets of type $\{r_{2t+1}\}$ or $\{r_{4t+2}\}$. Thus, $t(G)=3t+2=[(3n+5)/4]$.

Suppose that $n=4t+3$. Then $$N(G)=\{2t+3,2t+5,\ldots,4t+3,4t+4,4t+6,\ldots,8t+6\},$$ so
$|N(G)|=3t+3$ and $t(G)\geqslant6$. Since $t(2,G)\leqslant t(p,G)<4$, primes $p$ and $2$ cannot
lie in a coclique of maximal size. If $\eta(i)<2t+1$, then $r_i$ is adjacent to $r_{4t+4}$,
$r_{4t+6}$, and $r_{2t+3}$. Assume that $\eta(i)=2t+1$. If $r_i$ is adjacent to $r_j$ with $j\in
N(G)$, then $j=4t+4$. Since there are two distinct numbers $2t+1$ and $4t+2$ such that the value
of function $\eta$ of them is equal to $2t+1$, we have that $\Theta'(G)$ consists of one-element
sets of one of three types: $\{r_{4t+4}\}$, $\{r_{2t+1}\}$ or $\{r_{4t+2}\}$. Thus,
$M(G)=N(G)\setminus\{4t+4\}$, every $\theta(G)\in\Theta(G)$ is of type $\{r_i\mid
(n+1)/2<\eta(i)\leqslant n\}$, and $t(G)=3t+3=[(3n+5)/4]$.

It remains to consider the cases: $q=2$ and $n=4+k$, where $k=0,1,2,3$. All results (see
Table~\ref{ClassicTable}) are obtained by arguments similar to that in general case with respect
to the fact: $R_{4t+2}=R_6=\varnothing$, and can be easily verified by using
Proposition~\ref{adjbn} and \cite[Propositions 3.1, 4.3]{VasVd}.

\textbf{Case 2.} Let $G=D_n(q)$.

Suppose that $n=4$. If $q\neq2$ then $N(G)=\{3,6\}$, and if $q=2$ then $N(G)=\{3\}$. The set
$\{1,2,3,4,6\}$ includes $I(G)$, so $\xi=\{p\}\cup R_1\cup R_2\cup R_4$, where $\{p\}$, $R_2$ and $R_4$
are always nonempty. The prime $p$ and any prime $r_4$ are adjacent one to another, and are
non-adjacent to every $r_i$ with $i\in N(G)$. On the other hand, for $i\in\{1,2\}$ and
$j\in\{3,6\}$, primes $r_i$ and $r_j$ are adjacent. Therefore, $M(G)=N(G)$, $\theta(G)$ is of
type $\{r_3\}$ for $q=2$, and is of type $\{r_3,r_6\}$ otherwise. The set $\Theta'(G)$ consists of
one-element sets of type $\{p\}$, $\{r_2\}$, and $\{r_4\}$, if $q=2$, and sets of type $\{p\}$, and
$\{r_4\}$ otherwise. Thus, $t(G)=2$ for $q=2$, and $t(G)=3$ otherwise.

Let $n>4$. Now we consider four different cases subject to residue of $n$ modulo $4$. We write
$n=4t+k$, where $k=0,1,2,3$, and $t\geqslant1$. If $q=2$ we assume that $t>1$ to avoid exceptional
cases that arise because of $R_6=\varnothing$ for $q=2$.

Suppose that $n=4t>4$. Then
$$N(G)=\{2t+1,2t+3,\ldots,4t-1,4t+2,4t+4,\ldots,8t-2\},$$
so $|N(G)|=3t-1>4$. By \cite[Propositions 3.1, 4.4]{VasVd}, we have $t(2,G)\leqslant t(p,G)<4$, so
$p$ and $2$ cannot lie in any coclique of maximal size. By adjacency criterion, $r_{4t}$ is
non-adjacent to every $r_i$, where $i\in N(G)$. On the other hand, any prime $r_{4t-2}$ is adjacent
to $r_{4t}$, $r_{4t+2}$, any prime $r_{2t-1}$ is adjacent to $r_{4t}$, $r_{2t+1}$, and if
$\eta(i)<2t-1$, then $r_i$ is adjacent to at least three primes from every $\chi$. Therefore,
$M(G)=N(G)\cup\{n\}$, every $\theta(G)\in\Theta(G)$ is of type $\{r_i\mid
n/2\leqslant\eta(i)\leqslant n,i\neq2n\}$, $\Theta'(G)=\varnothing$, and $t(G)=3t=[(3n+1)/4]$.

Suppose that $n=4t+1$. Then
$$N(G)=\{2t+1,2t+3,\ldots,4t+1,4t+2,4t+4,\ldots,8t\},$$
so $|N(G)|=3t+1\geqslant4$. By \cite[Propositions 3.1, 4.4]{VasVd}, we have $t(2,G)\leqslant
t(p,G)<4$. Therefore, $p$ and $2$ cannot lie in any coclique of maximal size. If $\eta(i)<2t$, then
any prime $r_i$ is adjacent to $r_{4t+2}$, $r_{2t+1}$. Assume that $i=4t$, then $r_i$ is adjacent
to $r_{j}$, where $j\in N(G)$, if and only if $j=4t+2$. Thus, $M(G)=N(G)\setminus\{n+1\}$, every
$\theta(G)\in\Theta(G)$ is of type $\{r_i\mid n/2<\eta(i)\leqslant n,i\neq n+1,2n\}$, and
$\Theta'(G)$ consists of one-element sets of type $\{r_{4t}\}$ or $\{r_{4t+2}\}$. Therefore,
$t(G)=3t+1=[(3n+1)/4]$.

Suppose that $n=4t+2$. Then
$$N(G)=\{2t+3,2t+5,\ldots,4t+1,4t+4,4t+6,\ldots,8t+2\},$$
so $|N(G)|=3t\geqslant3$. Any primes $r_{2t+1}$ and $r_{4t+2}$ are adjacent one to another and are
non-adjacent to every $r_i$ with $i\in N(G)$. Hence $t(G)\geqslant4$. Since $t(2,G)\leqslant
t(p,G)<4$, primes $p$ and $2$ cannot lie in any coclique of maximal size. If $\eta(i)<n/2$, then
$r_i$ is adjacent to $r_{2t+1}$, $r_{4t+2}$, $r_{4t+4}$. Therefore, $N(G)=M(G)$, every
$\theta(G)\in\Theta(G)$ is of type $\{r_i\mid n/2<\eta(i)\leqslant n, i\neq2n\}$, and $\Theta'(G)$
consists of one-element sets of type $\{r_{2t+1}\}$ or $\{r_{4t+2}\}$. Thus,
$t(G)=3t+1=[(3n+1)/4]$.

Suppose that $n=4t+3$. Then $$N(G)=\{2t+3,2t+5,\ldots,4t+3,4t+4,4t+6,\ldots,8t+4\},$$ so
$|N(G)|=3t+2\geqslant5$. Since $t(2,G)\leqslant t(p,G)<4$, primes $p$ and $2$ cannot lie in a
coclique of maximal size. By adjacency criterion, $r_{2t+1}$ is non-adjacent to every $r_i$, where
$i\in N(G)$. On the other hand, if $\eta(i)<2t+1$ or $i=4t+2$, then $r_i$ is adjacent to
$r_{4t+4}$, $r_{2t+1}$. Therefore, $M(G)=N(G)\cup\{(n-1)/2\}$, every $\theta(G)\in\Theta(G)$ is of
type $\{r_i\mid (n-1)/2\leqslant\eta(i)\leqslant n,i\neq2n,n-1\}$, $\Theta'(G)=\varnothing$, and
$t(G)=3t+3=(3n+3)/4$. In case $n=4t+3$, any coclique of maximal size does not contain primes of
type $r_{4t+2}$, and so the group $D_7(2)$ is considered as well.

It remains to consider the cases: $q=2$ and $n=4+k$, where $k=1,2$. Both results (see
Table~\ref{ClassicTable}) are obtained by arguments similar to that in general case with respect to
the fact: $R_6=\varnothing$, and can be easily verified by using Proposition~\ref{adjdn} and
\cite[Propositions 3.1, 4.4]{VasVd}.

\textbf{Case 3.} Let $G={}^2D_n(q)$.

Suppose that $n=4$. If $q\neq2$ then $N(G)=\{3,6,8\}$, and if $q=2$ then $N(G)=\{3,8\}$. The set
$\{1,2,3,4,6,8\}$ includes $I(G)$, and so $\xi=\{p\}\cup R_1\cup R_2\cup R_4$, where $\{p\}$, $R_2$, and
$R_4$ are always nonempty. The prime $p$ and any prime $r_4$ are adjacent one to another, and are
non-adjacent to every $r_i$ with $i\in N(G)$. On the other hand, for $i\in\{1,2\}$ and
$j\in\{3,6\}$, primes $r_i$ and $r_j$ are adjacent. Therefore, $M(G)=N(G)$, $\theta(G)$ is of type
$\{r_3,r_8\}$ for $q=2$, and is of type $\{r_3,r_6,r_8\}$ otherwise. The set $\Theta'(G)$ consists
of one-element sets of type $\{p\}$, and $\{r_4\}$. Thus, $t(G)=3$ for $q=2$, and $t(G)=4$
otherwise.

Let $n>4$. Now we consider four different cases subject to residue of $n$ modulo $4$. We write
$n=4t+k$, where $k=0,1,2,3$, and $t\geqslant1$. If $q=2$ we assume that $t>1$ to avoid exceptional
cases that arise because of $R_6=\varnothing$ for $q=2$.

Suppose that $n=4t>4$. Then
$$N(G)=\{2t+1,2t+3,\ldots,4t-1,4t+2,4t+4,\ldots,8t\},$$
so $|N(G)|=3t>4$. By \cite[Propositions 3.1, 4.4]{VasVd}, we have $t(2,G)\leqslant
t(p,G)\leqslant4$, so $p$ and $2$ cannot lie in any coclique of maximal size. By adjacency
criterion, $r_{4t}$ is non-adjacent to every $r_i$, where $i\in N(G)$. On the other hand, any prime
$r_{2t-1}$ is adjacent to $r_{4t}$, $r_{4t+2}$, any prime $r_{4t-2}$ is adjacent to $r_{4t}$,
$r_{2t+1}$, and if $\eta(i)<2t-1$, then $r_i$ is adjacent to at least three primes from every
$\chi$. Therefore, $M(G)=N(G)\cup\{n\}$, every $\theta(G)\in\Theta(G)$ is of type $\{r_i\mid
n/2\leqslant\eta(i)\leqslant n\}$, $\Theta'(G)=\varnothing$, and $t(G)=3t+1=[(3n+4)/4]$.

Suppose that $n=4t+1$. Then
$$N(G)=\{2t+1,2t+3,\ldots,4t-1,4t+2,4t+4,\ldots,8t+2\},$$
so $|N(G)|=3t+1\geqslant4$. By \cite[Propositions 3.1, 4.4]{VasVd}, we have $t(2,G)\leqslant
t(p,G)<4$. Therefore, $p$ and $2$ cannot lie in any coclique of maximal size. If $\eta(i)<2t$, then
any prime $r_i$ is adjacent to $r_{4t+2}$, $r_{2t+1}$. Assume that $i=4t$, then $r_i$ is adjacent
to $r_{j}$, where $j\in N(G)$, if and only if $j=2t+1$. Thus, $M(G)=N(G)\setminus\{(n+1)/2\}$,
every $\theta(G)\in\Theta(G)$ is of type $\{r_i\mid n/2<\eta(i)\leqslant n,i\neq (n+1)/2,n\}$, and
$\Theta'(G)$ consists of one-element sets of type $\{r_{4t}\}$ or $\{r_{2t+1}\}$. Therefore,
$t(G)=3t+1=[(3n+4)/4]$.

Suppose that $n=4t+2$. Then
$$N(G)=\{2t+3,2t+5,\ldots,4t+1,4t+4,4t+6,\ldots,8t+4\},$$
so $|N(G)|=3t+1\geqslant4$. Any primes $r_{2t+1}$, $r_{4t}$ and $r_{4t+2}$ are adjacent one to
another and are non-adjacent to every $r_i$ with $i\in N(G)$. Hence $t(G)>4$. Since
$t(2,G)\leqslant t(p,G)\leqslant4$, primes $p$ and $2$ cannot lie in any coclique of maximal size.
If $\eta(i)<2t$, then $r_i$ is adjacent to $r_{2t+1}$, $r_{4t+2}$, $r_{4t}$, $r_{4t+4}$. Therefore,
$N(G)=M(G)$, every $\theta(G)\in\Theta(G)$ is of type $\{r_i\mid n/2<\eta(i)\leqslant n\}$, and
$\Theta'(G)$ consists of one-element sets of type $\{r_{2t+1}\}$, $\{r_{4t}\}$ or $\{r_{4t+2}\}$.
Thus, $t(G)=3t+2=[(3n+4)/4]$.

Suppose that $n=4t+3$. Then $$N(G)=\{2t+3,2t+5,\ldots,4t+1,4t+4,4t+6,\ldots,8t+6\},$$ so
$|N(G)|=3t+2\geqslant5$. Since $t(2,G)\leqslant t(p,G)<4$, primes $p$ and $2$ cannot lie in a
coclique of maximal size. By adjacency criterion, $r_{4t+2}$ is non-adjacent to every $r_i$, where
$i\in N(G)$. On the other hand, if $\eta(i)<2t+1$ or $i=2t+1$, then $r_i$ is adjacent to
$r_{4t+4}$, $r_{4t+2}$. Therefore, $M(G)=N(G)\cup\{n-1\}$, every $\theta(G)\in\Theta(G)$ is of type
$\{r_i\mid (n-1)/2\leqslant\eta(i)\leqslant n,i\neq n,(n-1)/2\}$, $\Theta'(G)=\varnothing$, and
$t(G)=3t+3=[(3n+4)/4]$.

It remains to consider the cases: $q=2$ and $n=4+k$, where $k=1,2,3$. All results (see
Table~\ref{ClassicTable}) are obtained by arguments similar to that in general case with respect
to the fact: $R_{4t+2}=R_6=\varnothing$, and can be easily verified by using
Proposition~\ref{adjdn} and \cite[Propositions 3.1, 4.4]{VasVd}.
\end{proof}

\begin{prop}\label{cocliqueexcept} If $G$ is an finite simple exceptional group of Lie type
over a field of characteristic~$p$, then $t(G)$, and the sets $\Theta(G)$, $\Theta'(G)$ are listed
in Table~{\em\ref{ExceptTable}}.
\end{prop}

\begin{proof}
We consider all types of  exceptional groups of Lie type separately. Following \cite{ZavGraph}, we
use the compact form of the prime graph $GK(G)$. By the compact form we mean a graph whose vertices
are labeled with marks $R_i$. A vertex labeled $R_i$ represents the clique of $GK(G)$ such that
every vertex in this clique labeled by a prime from~$R_i$. An edge connecting $R_i$ and $R_j$
represents the set of edges of $GK(G)$ that connect each vertex in $R_i$ with each vertex in $R_j$.
If an edge occurs under some condition, we draw such edge by a dotted line and write corresponding
occurence condition. The technical tools for determining the compact form of the prime graph
$GK(G)$ for an exceptional group of Lie type $G$ are Propositions \ref{adjexcept} and
\ref{adjsuzree}, and also \cite[Propositions~3.2, 3.3, and~4.5]{VasVd}. Notice that the compact
form of $GK(G)$ can be considered as a graphical form of the adjacency criterion in~$GK(G)$.

\begin{center}
The compact form for $GK(G_2(q))$
\end{center}

%\vspace{1\baselineskip}
\setlength{\unitlength}{4144sp}%
\begingroup\makeatletter\ifx\SetFigFontNFSS\undefined%
\gdef\SetFigFontNFSS#1#2#3#4#5{%
  \reset@font\fontsize{#1}{#2pt}%
  \fontfamily{#3}\fontseries{#4}\fontshape{#5}%
  \selectfont}%
\fi\endgroup%
\begin{picture}(2724,1374)(2014,-2548)
\thinlines
{\color[rgb]{0,0,0}\put(2701,-2536){\line( 1, 0){1350}}
\put(4051,-2536){\line(-1, 1){675}}
\put(3376,-1861){\line(-1,-1){675}}
}%
{\color[rgb]{0,0,0}\put(2026,-1861){\line( 1,-1){675}}
\color[rgb]{0,0,0}\put(2701,-2536){\line(-2,1){180}}
\color[rgb]{0,0,0}\put(2701,-2536){\line(-1,2){90}}
}%
{\color[rgb]{0,0,0}\put(4726,-1861){\line(-1,-1){675}}
\color[rgb]{0,0,0}\put(4051,-2536){\line(2,1){180}}
\color[rgb]{0,0,0}\put(4051,-2536){\line(1,2){90}}
}%
{\color[rgb]{0,0,0}\multiput(2026,-1861)(0.00000,122.72727){6}{\line( 0, 1){ 61.364}}
}%
{\color[rgb]{0,0,0}\multiput(4726,-1861)(0.00000,122.72727){6}{\line( 0, 1){ 61.364}}
}%
\put(2026,-1861){\circle*{60}}\put(2126,-1801){\makebox(0,0){$3$}}
\put(4726,-1861){\circle*{60}}\put(4636,-1801){\makebox(0,0){$3$}}
\put(4726,-1181){\circle*{60}}\put(4566,-1181){\makebox(0,0){$R_6$}}
\put(4051,-2536){\circle*{60}} \put(4011,-2306){\makebox(0,0){$R_2$}}
\put(2701,-2536){\circle*{60}} \put(2751,-2306){\makebox(0,0){$R_1$}}
\put(3376,-1861){\circle*{60}} \put(3376,-1731){\makebox(0,0){$p$}}
\put(2026,-1181){\circle*{60}} \put(2186,-1181){\makebox(0,0){$R_3$}}
\end{picture}%
\vspace{1\baselineskip}

Let $G=G_2(q)$. In the compact form for $GK(G_2(q))$ the vector from $3$ to $R_1$ (resp. $R_2$) and
the dotted edge $(3,R_{3})$ (resp.$(3,R_{6})$)  mean that $R_1$ (resp. $R_2$) and $R_{3}$ (resp.
$R_6$) are not connected, but if $3\in R_1$, i.~e., $q\equiv 1\pmod3$ (resp. $3\in R_2$, i.~e.,
$q\equiv-1\pmod 3$), then there exists an edge between $3$ and $R_{3}$ (resp. $R_6$). If
$R_1=\varnothing$, then one need to remove vertex $R_1$ with all corresponding edges. From the
compact form of $GK(G)$ it is evident, that $\Theta(G)=\{\{r_3,r_6\}\mid r_i\in R_i\}$, while
$\Theta'(G)=\{\{p\}, \{r_1\},\{r_2\}\mid r_i\in R_i\setminus\{3\}\}$.

 \begin{center}
The compact form for $GK(F_4(q))$
\end{center}
\begin{tabular}{ccc}
\setlength{\unitlength}{3108sp}%
\begingroup\makeatletter\ifx\SetFigFontNFSS\undefined%
\gdef\SetFigFontNFSS#1#2#3#4#5{%
  \reset@font\fontsize{#1}{#2pt}%
  \fontfamily{#3}\fontseries{#4}\fontshape{#5}%
  \selectfont}%
\fi\endgroup%
\begin{picture}(2724,2724)(3139,-3898)
\thinlines
{\color[rgb]{0,0,0}\put(5851,-1861){\line(-1, 1){675}}
}%
\put(5851,-1861){\circle*{90}} \put(6281,-1861){\makebox(0,0){$2(\not=p)$}}
\put(3826,-1186){\circle*{90}} \put(3576,-1186){\makebox(0,0){$R_4$}}
\put(4535,-1186){\circle*{90}} \put(4735,-1186){\makebox(0,0){$R_{12}$}}
{\color[rgb]{0,0,0}\put(5851,-1861){\line(-3, 1){2025}}
}%
\put(5176,-1186){\circle*{90}} \put(5376,-1186){\makebox(0,0){$R_8$}}
{\color[rgb]{0,0,0}\put(5851,-1861){\line(-1, 0){2700}}
}%
{\color[rgb]{0,0,0}\put(5851,-1861){\line(-1,-1){2025}}
}%
{\color[rgb]{0,0,0}\put(5851,-1861){\line(-1,-3){675}}
}%
{\color[rgb]{0,0,0}\put(5851,-1861){\line( 0,-1){1350}}
}%
{\color[rgb]{0,0,0}\put(3151,-1861){\line( 1, 1){675}}
}%
\put(3151,-1861){\circle*{90}} \put(3001,-1861){\makebox(0,0){$p$}}
{\color[rgb]{0,0,0}\put(3151,-1861){\line( 0,-1){1350}}
}%
\put(3151,-3211){\circle*{90}} \put(2951,-3211){\makebox(0,0){$R_1$}}
{\color[rgb]{0,0,0}\put(3151,-1861){\line( 1,-3){675}}
}%
{\color[rgb]{0,0,0}\put(3151,-1861){\line( 1,-1){2025}}
}%
{\color[rgb]{0,0,0}\put(3151,-3211){\line( 1,-1){675}}
}%
{\color[rgb]{0,0,0}\put(3151,-3211){\line( 3,-1){2025}}
}%
\put(5851,-3211){\circle*{90}} \put(6051,-3211){\makebox(0,0){$R_2$}}
{\color[rgb]{0,0,0}\put(5851,-3211){\line(-1,-1){675}}
}%
{\color[rgb]{0,0,0}\put(5851,-3211){\line(-3,-1){2025}}
}%
{\color[rgb]{0,0,0}\put(3826,-1186){\line(-1,-3){675}}
}%
{\color[rgb]{0,0,0}\put(3151,-3211){\line( 1, 0){2700}}
\put(5851,-3211){\line(-1, 1){2025}}
}%
{\color[rgb]{0,0,0}\put(3151,-3211){\line( 2, 1){2700}}
}%
{\color[rgb]{0,0,0}\put(5851,-3211){\line(-2, 1){2700}}
}%
\put(3826,-3886){\circle*{90}} \put(3626,-3886){\makebox(0,0){$R_3$}}
\put(5176,-3886){\circle*{90}} \put(5376,-3886){\makebox(0,0){$R_6$}}
\end{picture}&\ \ \ \ \ \ \ \ \ \ \ \ \ \ \ \ \  &
\setlength{\unitlength}{3108sp}%
\begingroup\makeatletter\ifx\SetFigFontNFSS\undefined%
\gdef\SetFigFontNFSS#1#2#3#4#5{%
  \reset@font\fontsize{#1}{#2pt}%
  \fontfamily{#3}\fontseries{#4}\fontshape{#5}%
  \selectfont}%
\fi\endgroup%
\begin{picture}(2724,2724)(3139,-3898)
\thinlines
%{\color[rgb]{0,0,0}\put(5851,-1861){\line(-1, 1){675}}
%}%
%\put(5851,-1861){\circle*{90}} \put(6281,-1861){\makebox(0,0){$2(\not=p)$}}
\put(3826,-1186){\circle*{90}} \put(3576,-1186){\makebox(0,0){$R_4$}}
\put(4535,-1186){\circle*{90}} \put(4735,-1186){\makebox(0,0){$R_{12}$}}
%{\color[rgb]{0,0,0}\put(5851,-1861){\line(-3, 1){2025}}
%}%
\put(5176,-1186){\circle*{90}} \put(5376,-1186){\makebox(0,0){$R_8$}}
%{\color[rgb]{0,0,0}\put(5851,-1861){\line(-1, 0){2700}}
%}%
%{\color[rgb]{0,0,0}\put(5851,-1861){\line(-1,-1){2025}}
%}%
%{\color[rgb]{0,0,0}\put(5851,-1861){\line(-1,-3){675}}
%}%
%{\color[rgb]{0,0,0}\put(5851,-1861){\line( 0,-1){1350}}
%}%
{\color[rgb]{0,0,0}\put(3151,-1861){\line( 1, 1){675}}
}%
\put(3151,-1861){\circle*{90}} \put(2721,-1861){\makebox(0,0){$(2=)p$}}
{\color[rgb]{0,0,0}\put(3151,-1861){\line( 0,-1){1350}}
}%
\put(3151,-3211){\circle*{90}} \put(2951,-3211){\makebox(0,0){$R_1$}}
{\color[rgb]{0,0,0}\put(3151,-1861){\line( 1,-3){675}}
}%
{\color[rgb]{0,0,0}\put(3151,-1861){\line( 1,-1){2025}}
}%
{\color[rgb]{0,0,0}\put(3151,-3211){\line( 1,-1){675}}
}%
{\color[rgb]{0,0,0}\put(3151,-3211){\line( 3,-1){2025}}
}%
\put(5851,-3211){\circle*{90}} \put(6051,-3211){\makebox(0,0){$R_2$}}
{\color[rgb]{0,0,0}\put(5851,-3211){\line(-1,-1){675}}
}%
{\color[rgb]{0,0,0}\put(5851,-3211){\line(-3,-1){2025}}
}%
{\color[rgb]{0,0,0}\put(3826,-1186){\line(-1,-3){675}}
}%
{\color[rgb]{0,0,0}\put(3151,-3211){\line( 1, 0){2700}}
\put(5851,-3211){\line(-1, 1){2025}}
}%
%{\color[rgb]{0,0,0}\put(3151,-3211){\line( 2, 1){2700}}
%}%
{\color[rgb]{0,0,0}\put(5851,-3211){\line(-2, 1){2700}}
}%
\put(3826,-3886){\circle*{90}} \put(3626,-3886){\makebox(0,0){$R_3$}}
\put(5176,-3886){\circle*{90}} \put(5376,-3886){\makebox(0,0){$R_6$}}
\end{picture}\\
\end{tabular}
\vspace{1\baselineskip}

Let $G=F_4(q)$. It is evident from the compact form for $GK(F_4(q))$
that $\{2, p,R_1,R_2,R_3\}$ is a clique, while the remaining
vertices in the compact form are pairwise non-adjacent. Since $R_3$ in non-adjacent to $R_4,R_6,R_8,R_{12}$ and the remaining vertices from
the set
$\{2,p, R_1,R_2\}$ are adjacent to at least two vertices from the set  $\{R_4,R_6,R_8,R_{12}\}$, we obtain that
 $\Theta(G)=\{\{r_3,r_4,r_6,r_8,r_{12}\}\mid r_i\in R_i\}$ if $R_6\not=\varnothing$ and
$\Theta(G)=\{\{r_3,r_4,r_8,r_{12}\}\mid r_i\in R_i\}$ if $R_6=\varnothing$ (i.~e., if $q=2$). In
both cases~${\Theta'(G)=\varnothing}$.
%\vspace{1\baselineskip}

\begin{center}
The compact form for $GK(E_6^\varepsilon(q))$
\end{center}

\vspace{1\baselineskip}

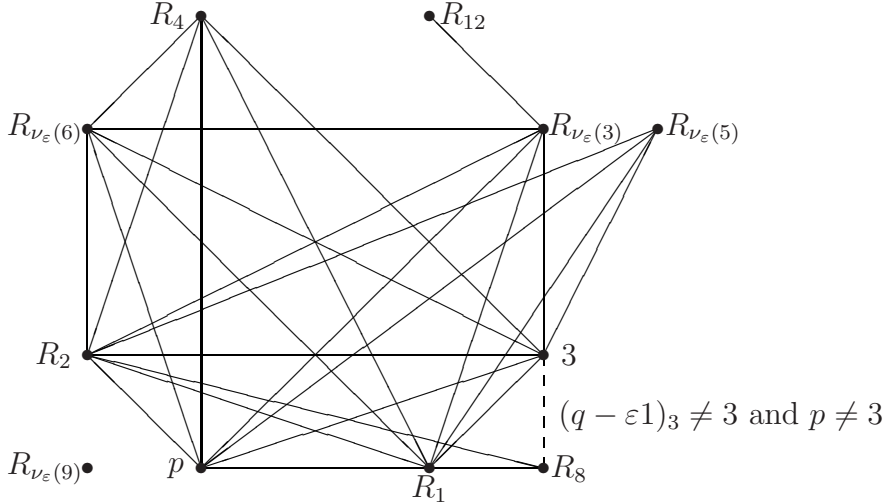
\setlength{\unitlength}{4144sp}%
\begingroup\makeatletter\ifx\SetFigFontNFSS\undefined%
\gdef\SetFigFontNFSS#1#2#3#4#5{%
  \reset@font\fontsize{#1}{#2pt}%
  \fontfamily{#3}\fontseries{#4}\fontshape{#5}%
  \selectfont}%
\fi\endgroup%
\begin{picture}(3399,2724)(2464,-3898)
\thinlines
{\color[rgb]{0,0,0}\put(5851,-3211){\line(-1, 1){2025}}
}%
{\color[rgb]{0,0,0}\put(5851,-3211){\line(-3,-1){2025}}
}%
{\color[rgb]{0,0,0}\put(5851,-3211){\line( 1, 2){675}}
}%
{\color[rgb]{0,0,0}\put(5851,-3211){\line(-2, 1){2700}}
}%
{\color[rgb]{0,0,0}\put(5851,-3211){\line(-1,-1){675}}
}%
{\color[rgb]{0,0,0}\put(5851,-3211){\line( 0, 1){1350}}
}%
{\color[rgb]{0,0,0}\put(6526,-1861){\line(-4,-3){2700}}
}%
{\color[rgb]{0,0,0}\put(6526,-1861){\line(-5,-2){3375}}
}%
{\color[rgb]{0,0,0}\put(6526,-1861){\line(-2,-3){1350}}
}%
{\color[rgb]{0,0,0}\put(3151,-3211){\line( 0, 1){1350}}
}%
{\color[rgb]{0,0,0}\put(3151,-3211){\line( 3,-1){2025}}
}%
{\color[rgb]{0,0,0}\put(5851,-1861){\line(-1,-3){675}}
}%
{\color[rgb]{0,0,0}\put(5851,-1861){\line(-1, 0){2700}}
}%
{\color[rgb]{0,0,0}\put(3826,-3886){\line(-1, 1){675}}
}%
{\color[rgb]{0,0,0}\put(3151,-3211){\line( 2, 1){2700}}
\put(5851,-1861){\line(-1,-1){2025}}
}%
{\color[rgb]{0,0,0}\put(3151,-3211){\line( 1, 0){2700}}
}%
{\color[rgb]{0,0,0}\put(3151,-1861){\line( 1,-1){2025}}
}%
{\color[rgb]{0,0,0}\put(5176,-3886){\line(-1, 0){1350}}
\put(3826,-3886){\line( 0, 1){2700}}
\put(3826,-1186){\line(-1,-3){675}}
}%
{\color[rgb]{0,0,0}\put(3826,-1186){\line( 1,-2){1350}}
}%
{\color[rgb]{0,0,0}\put(3826,-1186){\line(-1,-1){675}}
}%
{\color[rgb]{0,0,0}\multiput(5851,-3211)(0.00000,-122.72727){6}{\line( 0,-1){ 61.364}}
}%
{\color[rgb]{0,0,0}\put(3151,-3211){\line( 4,-1){2700}}
\put(5851,-3886){\line(-1, 0){675}}
}%
{\color[rgb]{0,0,0}\put(5851,-1861){\line(-1, 1){675}}
}%
{\color[rgb]{0,0,0}\put(3826,-3886){\line(-1, 3){675}}
}%
\put(3826,-1186){\circle*{60}} \put(3626,-1186){\makebox(0,0){$R_4$}}
\put(3826,-3886){\circle*{60}} \put(3676,-3886){\makebox(0,0){$p$}}
\put(5176,-1186){\circle*{60}} \put(5376,-1186){\makebox(0,0){$R_{12}$}}
\put(3151,-1861){\circle*{60}} \put(2901,-1861){\makebox(0,0){$R_{\nu_\varepsilon(6)}$}}
%\put(5851,-1186){\circle*{60}} \put(6001,-1186){\makebox(0,0){$R_8$}}
\put(6901,-3576){\makebox(0,0){$(q-\varepsilon1)_3\not=3$ and $p\not=3$}}
\put(5851,-3211){\circle*{60}} \put(6001,-3211){\makebox(0,0){$3$}}
\put(3151,-3211){\circle*{60}} \put(2951,-3211){\makebox(0,0){$R_2$}}
\put(5851,-1861){\circle*{60}} \put(6101,-1861){\makebox(0,0){$R_{\nu_\varepsilon(3)}$}}
\put(6526,-1861){\circle*{60}} \put(6796,-1861){\makebox(0,0){$R_{\nu_\varepsilon(5)}$}}
\put(5176,-3886){\circle*{60}} \put(5176,-4006){\makebox(0,0){$R_1$}}
\put(5851,-3886){\circle*{60}} \put(6001,-3886){\makebox(0,0){$R_{8}$}}
\put(3151,-3886){\circle*{60}} \put(2901,-3886){\makebox(0,0){$R_{\nu_\varepsilon(9)}$}}
\end{picture}%
\vspace{1\baselineskip}

Let $G=E_6^\varepsilon(q)$. In the compact form for $GK(E_6(q))$ the set $\{3, p,
R_1,R_2,R_{\nu_\varepsilon(3)},R_{\nu_\varepsilon(6)}\}$ forms a clique, while the remaining
vertices are pairwise non-adjacent. Moreover, $R_{\nu_\varepsilon(3)}$ and $R_{\nu_\varepsilon(6)}$
are the only vertices from $\{3, p, R_1,R_2,R_{\nu_\varepsilon(3)},R_{\nu_\varepsilon(6)}\}$, which
are adjacent to precisely one of the remaining vertices (namely, $R_{\nu_\varepsilon(3)}$ is
adjacent to $R_{12}$, and $R_{\nu_\varepsilon(6)}$ is adjacent to $R_{4}$). Thus
$\Theta(G)=\{\{r_{\nu_\varepsilon(5)},r_8,r_{\nu_\varepsilon(9)}\}\mid r_i\in R_i\}$ and
$\Theta'(G)=\{\{r_4,r_{\nu_\varepsilon(3)}\},\{r_{\nu_\varepsilon(6)},r_{12}\}, \{r_4,r_{12}\}\mid
r_i\in R_i\}$. Since $R_6=\varnothing$ for $q=2$, we obtain exceptions mentioned in
Table~\ref{ExceptTable}.

\begin{center}
The compact form for $GK(E_7(q))$
\end{center}
\setlength{\unitlength}{3108sp}%
\begingroup\makeatletter\ifx\SetFigFontNFSS\undefined%
\gdef\SetFigFontNFSS#1#2#3#4#5{%
  \reset@font\fontsize{#1}{#2pt}%
  \fontfamily{#3}\fontseries{#4}\fontshape{#5}%
  \selectfont}%
\fi\endgroup%
\begin{picture}(5424,3624)(889,-4123)
\thinlines
{\color[rgb]{0,0,0}\put(5401,-1411){\line( 0,-1){1800}}
}%
{\color[rgb]{0,0,0}\put(5401,-1411){\line(-1,-3){900}}
}%
{\color[rgb]{0,0,0}\put(5401,-1411){\line(-1, 0){3600}}
}%
{\color[rgb]{0,0,0}\put(1801,-3211){\line( 2, 1){3600}}
}%
{\color[rgb]{0,0,0}\put(5401,-1411){\line(-1,-1){2700}}
}%
{\color[rgb]{0,0,0}\put(1801,-3211){\line( 1, 0){3600}}
}%
{\color[rgb]{0,0,0}\put(5401,-3211){\line(-2, 1){3600}}
}%
{\color[rgb]{0,0,0}\put(5401,-3211){\line( 1, 2){900}}
}%
{\color[rgb]{0,0,0}\put(5401,-3211){\line(-3,-1){2700}}
}%
{\color[rgb]{0,0,0}\put(5401,-3211){\line(-1,-1){900}}
}%
{\color[rgb]{0,0,0}\put(4501,-511){\line( 0,-1){3600}}
}%
{\color[rgb]{0,0,0}\put(4501,-511){\line(-1,-2){1800}}
}%
{\color[rgb]{0,0,0}\put(1801,-3211){\line( 1, 1){2700}}
}%
{\color[rgb]{0,0,0}\put(1801,-3211){\line( 3,-1){2700}}
}%
{\color[rgb]{0,0,0}\put(1801,-3211){\line( 5, 2){4500}}
}%
{\color[rgb]{0,0,0}\put(1801,-3211){\line( 1,-1){900}}
}%
{\color[rgb]{0,0,0}\put(1801,-1411){\line( 0,-1){1800}}
}%
{\color[rgb]{0,0,0}\put(4501,-4111){\line(-1, 1){2700}}
}%
{\color[rgb]{0,0,0}\put(2701,-4111){\line(-1, 3){900}}
}%
{\color[rgb]{0,0,0}\put(4501,-4111){\line(-1, 0){1800}}
}%
{\color[rgb]{0,0,0}\put(4501,-4111){\line( 2, 3){1800}}
}%
{\color[rgb]{0,0,0}\put(6301,-1411){\line(-4,-3){3600}}
}%
{\color[rgb]{0,0,0}\put(901,-1411){\line( 1,-2){900}}
}%
{\color[rgb]{0,0,0}\put(901,-1411){\line( 2,-3){1800}}
}%
{\color[rgb]{0,0,0}\put(901,-1411){\line( 4,-3){3600}}
}%
{\color[rgb]{0,0,0}\put(2701,-511){\line(-1,-3){900}}
}%
{\color[rgb]{0,0,0}\put(2701,-511){\line(-1,-1){900}}
}%
{\color[rgb]{0,0,0}\put(2701,-511){\line( 1,-1){2700}}
}%
{\color[rgb]{0,0,0}\put(2701,-511){\line( 0,-1){3600}}
}%
{\color[rgb]{0,0,0}\put(2701,-511){\line( 1,-2){1800}}
}%
{\color[rgb]{0,0,0}\put(4501,-4111){\line( 1, 0){1800}}
}%
{\color[rgb]{0,0,0}\put(4501,-4111){\line( 2, 1){1800}}
}%
{\color[rgb]{0,0,0}\put(2701,-4111){\line(-1, 0){1800}}
}%
{\color[rgb]{0,0,0}\put(2701,-4111){\line(-2, 1){1800}}
}%
{\color[rgb]{0,0,0}\put(4501,-511){\line( 1,-1){900}}
}%
{\color[rgb]{0,0,0}\put(901,-1411){\line( 1, 0){900}}
}%
\put(1801,-1411){\circle*{60}} \put(1621,-1301){\makebox(0,0){$R_6$}}
\put(4501,-511){\circle*{60}} \put(4701,-511){\makebox(0,0){$R_8$}}
\put(5401,-3211){\circle*{60}} \put(5561,-3211){\makebox(0,0){$R_3$}}
\put(2701,-4111){\circle*{60}} \put(2701,-4261){\makebox(0,0){$R_1$}}
\put(5401,-1411){\circle*{60}} \put(5451,-1131){\makebox(0,0){$R_4$}}
\put(1801,-3211){\circle*{60}} \put(1801,-3411){\makebox(0,0){$p$}}
\put(4501,-4111){\circle*{60}} \put(4501,-4261){\makebox(0,0){$R_2$}}
\put(901,-4111){\circle*{60}} \put(721,-4111){\makebox(0,0){$R_9$}}
\put(6301,-4111){\circle*{60}} \put(6501,-4111){\makebox(0,0){$R_{14}$}}
\put(6301,-3211){\circle*{60}} \put(6501,-3211){\makebox(0,0){$R_{18}$}}
\put(2701,-511){\circle*{60}} \put(2401,-511){\makebox(0,0){$R_{12}$}}
\put(901,-1411){\circle*{60}} \put(681,-1411){\makebox(0,0){$R_{10}$}}
\put(6301,-1411){\circle*{60}} \put(6501,-1411){\makebox(0,0){$R_{5}$}}
\put(901,-3211){\circle*{60}} \put(721,-3211){\makebox(0,0){$R_{7}$}}
\end{picture}%
\vspace{1\baselineskip}

Let $G=E_7(q)$. In the compact form for $GK(E_7(q))$ the set $\{p, R_1,R_2,R_3,R_4,R_6\}$ forms a
clique, while the remaining vertices are pairwise non-adjacent. Moreover, $R_4$ is the only
vertices from $\{p, R_1,R_2,R_3,R_4,R_6\}$, which are adjacent to precisely one of the remaining
vertices (namely, $R_4$ is adjacent to $R_8$). Thus
$\Theta(G)=\{\{r_5,r_7,r_9,r_{10},r_{12},r_{14},r_{18}\}\mid r_i\in R_i\}$ and
$\Theta'(G)=\{\{r_4\},\{r_8\}\mid r_i\in R_i\}$.

%\newpage

\begin{center}
The compact form for $GK(E_8(q))$
\end{center}
\setlength{\unitlength}{4144sp}%
\begingroup\makeatletter\ifx\SetFigFontNFSS\undefined%
\gdef\SetFigFontNFSS#1#2#3#4#5{%
  \reset@font\fontsize{#1}{#2pt}%
  \fontfamily{#3}\fontseries{#4}\fontshape{#5}%
  \selectfont}%
\fi\endgroup%
\begin{picture}(4749,4074)(2464,-4573)
\thinlines
{\color[rgb]{0,0,0}\put(5851,-1861){\line(-1, 2){675}}
}%
{\color[rgb]{0,0,0}\put(5851,-1861){\line(-3, 1){2025}}
}%
{\color[rgb]{0,0,0}\put(5851,-1861){\line(-1, 0){2700}}
}%
{\color[rgb]{0,0,0}\put(5851,-1861){\line(-1,-1){2025}}
}%
{\color[rgb]{0,0,0}\put(5851,-1861){\line( 0,-1){1350}}
}%
{\color[rgb]{0,0,0}\put(3151,-1861){\line( 1, 1){675}}
}%
{\color[rgb]{0,0,0}\put(3151,-1861){\line( 0,-1){1350}}
}%
{\color[rgb]{0,0,0}\put(3151,-3211){\line( 1,-1){675}}
}%
{\color[rgb]{0,0,0}\put(5851,-3211){\line(-3,-1){2025}}
}%
{\color[rgb]{0,0,0}\put(3826,-1186){\line(-1,-3){675}}
}%
{\color[rgb]{0,0,0}\put(3151,-3211){\line( 1, 0){2700}}
\put(5851,-3211){\line(-1, 1){2025}}
}%
{\color[rgb]{0,0,0}\put(3151,-3211){\line( 2, 1){2700}}
}%
{\color[rgb]{0,0,0}\put(5176,-511){\line(-2,-5){1350}}
}%
{\color[rgb]{0,0,0}\put(3826,-1186){\line( 0,-1){2700}}
}%
{\color[rgb]{0,0,0}\put(3151,-1861){\line( 1,-3){675}}
}%
{\color[rgb]{0,0,0}\put(3151,-1861){\line( 2,-1){2700}}
}%
{\color[rgb]{0,0,0}\put(3151,-1861){\line( 3, 2){2025}}
}%
{\color[rgb]{0,0,0}\put(3151,-3211){\line( 3, 4){2025}}
}%
{\color[rgb]{0,0,0}\put(5851,-3211){\line(-1, 4){675}}
}%
{\color[rgb]{0,0,0}\put(6526,-4561){\line(-1, 2){675}}
}%
{\color[rgb]{0,0,0}\put(6526,-4561){\line(-4, 5){2700}}
}%
{\color[rgb]{0,0,0}\put(6526,-4561){\line(-5, 4){3375}}
}%
{\color[rgb]{0,0,0}\put(6526,-4561){\line(-5, 2){3375}}
}%
{\color[rgb]{0,0,0}\put(6526,-4561){\line(-4, 1){2700}}
}%
{\color[rgb]{0,0,0}\put(2476,-511){\line( 1,-2){675}}
}%
{\color[rgb]{0,0,0}\put(2476,-511){\line( 1,-4){675}}
}%
{\color[rgb]{0,0,0}\put(2476,-511){\line( 5,-4){3375}}
}%
{\color[rgb]{0,0,0}\put(2476,-511){\line( 5,-2){3375}}
}%
{\color[rgb]{0,0,0}\put(6526,-511){\line(-4,-1){2700}}
}%
{\color[rgb]{0,0,0}\put(6526,-511){\line(-5,-2){3375}}
}%
{\color[rgb]{0,0,0}\put(6526,-511){\line(-5,-4){3375}}
}%
{\color[rgb]{0,0,0}\put(6526,-511){\line(-4,-5){2700}}
}%
{\color[rgb]{0,0,0}\put(5851,-1861){\line( 1, 2){675}}
}%
{\color[rgb]{0,0,0}\put(2476,-4561){\line( 2, 1){1350}}
}%
{\color[rgb]{0,0,0}\put(2476,-1861){\line( 1, 0){675}}
}%
{\color[rgb]{0,0,0}\put(2476,-1861){\line( 5,-2){3375}}
}%
{\color[rgb]{0,0,0}\put(2476,-1861){\line( 1,-2){675}}
}%
{\color[rgb]{0,0,0}\put(3151,-1861){\line(-1,-1){675}}
}%
{\color[rgb]{0,0,0}\put(2476,-2536){\line( 5,-1){3375}}
}%
{\color[rgb]{0,0,0}\put(2476,-2536){\line( 1,-1){675}}
}%
{\color[rgb]{0,0,0}\put(2476,-4561){\line( 5, 2){3375}}
}%
{\color[rgb]{0,0,0}\put(2476,-4561){\line( 1, 2){675}}
}%
{\color[rgb]{0,0,0}\put(2476,-4561){\line( 1, 4){675}}
}%
{\color[rgb]{0,0,0}\put(5851,-1861){\line( 1,-4){675}}
}%
{\color[rgb]{0,0,0}\put(5851,-3211){\line( 1, 4){675}}
}%
{\color[rgb]{0,0,0}\put(3826,-1186){\line(-2, 1){1350}}
}%
{\color[rgb]{0,0,0}\put(3826,-511){\line( 0,-1){675}}
}%
{\color[rgb]{0,0,0}\put(3826,-511){\line(-1,-4){675}}
}%
{\color[rgb]{0,0,0}\put(3826,-511){\line(-1,-2){675}}
}%
{\color[rgb]{0,0,0}\put(3826,-511){\line( 3,-4){2025}}
}%
{\color[rgb]{0,0,0}\put(6526,-1861){\line(-1, 0){675}}
}%
{\color[rgb]{0,0,0}\multiput(6526,-1861)(54.00000,0.00000){13}{\line( 1, 0){ 27.000}}
}%
{\color[rgb]{0,0,0}\multiput(6526,-1861)(54.00000,0.00000){13}{\line( 1, 0){ 27.000}}
\color[rgb]{0,0,0}\put(5851,-1861){\line(3, 1){180}}
\color[rgb]{0,0,0}\put(5851,-1861){\line(3, -1){180}}
}%
\put(3826,-1186){\circle*{60}} \put(3986,-1026){\makebox(0,0){$R_6$}}
\put(5176,-511){\circle*{60}} \put(5326,-511){\makebox(0,0){$R_5$}}
%\put(6851,-1711){\makebox(0,0){$q^2\equiv -1\pmod5$}}
\put(5851,-1861){\circle*{60}}\put(5841,-1591){\makebox(0,0){$R_4$}}
\put(6526,-1861){\circle*{60}}\put(6526,-1971){\makebox(0,0){$5$}}
\put(7201,-1861){\circle*{60}} \put(7201,-1971){\makebox(0,0){$R_{20}$}}
\put(3151,-1861){\circle*{60}} \put(3151,-1611){\makebox(0,0){$R_1$}}
\put(3151,-3211){\circle*{60}} \put(3021,-3181){\makebox(0,0){$R_2$}}
\put(5851,-3211){\circle*{60}} \put(5941,-3181){\makebox(0,0){$p$}}
\put(3826,-3886){\circle*{60}} \put(3656,-3836){\makebox(0,0){$R_3$}}
\put(7201,-3886){\circle*{60}} \put(7201,-3996){\makebox(0,0){$R_{24}$}}
\put(7201,-3211){\circle*{60}} \put(7201,-3321){\makebox(0,0){$R_{15}$}}
\put(7201,-2536){\circle*{60}} \put(7201,-2646){\makebox(0,0){$R_{30}$}}
\put(2476,-2536){\circle*{60}} \put(2321,-2536){\makebox(0,0){$R_{14}$}}
\put(2476,-1861){\circle*{60}} \put(2321,-1861){\makebox(0,0){$R_{7}$}}
\put(6526,-4561){\circle*{60}} \put(6676,-4561){\makebox(0,0){$R_{12}$}}
\put(6526,-511){\circle*{60}} \put(6626,-511){\makebox(0,0){$R_{8}$}}
\put(2476,-511){\circle*{60}} \put(2331,-511){\makebox(0,0){$R_{10}$}}
\put(3826,-511){\circle*{60}} \put(3986,-511){\makebox(0,0){$R_{18}$}}
\put(2476,-4561){\circle*{60}} \put(2326,-4561){\makebox(0,0){$R_{9}$}}
\end{picture}%
\vspace{1\baselineskip}

Let $G=E_8(q)$. In the compact form  for $GK(E_8(q))$, the vector from $5$ to $R_4$ and the dotted
edge $(5,R_{20})$ mean that $R_4$ and $R_{20}$ are not connected, but if $5\in R_4$ (i.~e.,
$q^2\equiv -1\pmod5$), then there exists an edge between $5$ and $R_{20}$. Now
$\{p,R_1,R_2,R_3,R_4,R_6\}$ forms a clique, while the remaining vertices are pairwise non-adjacent.
Notice that each vertex from the clique $\{p,R_1,R_2,R_3,R_4,R_6\}$ is adjacent to at least two
vertices from the set of remaining vertices. So
$$\Theta(G)=\{\{r_5,r_7,r_8,r_9,r_{10}, r_{12}, r_{14},r_{15},r_{18},r_{20}, r_{24}, r_{30}\}\mid r_i\in R_i\}$$ and
$\Theta'(G)=\varnothing$.

\begin{center}
The compact form for $GK({}^3D_4(q))$
\end{center}
\setlength{\unitlength}{4144sp}%
\begingroup\makeatletter\ifx\SetFigFontNFSS\undefined%
\gdef\SetFigFontNFSS#1#2#3#4#5{%
  \reset@font\fontsize{#1}{#2pt}%
  \fontfamily{#3}\fontseries{#4}\fontshape{#5}%
  \selectfont}%
\fi\endgroup%
\begin{picture}(1374,1374)(2014,-1423)
\thinlines
{\color[rgb]{0,0,0}\put(2701,-61){\line(-1,-1){675}}
\put(2026,-736){\line( 0,-1){675}}
\put(2026,-1411){\line( 2, 1){1350}}
\put(3376,-736){\line(-1, 1){675}}
\put(2701,-61){\line(-1,-2){675}}
}%
\put(2701,-61){\circle*{60}} \put(2851,-61){\makebox(0,0){$p$}}
\put(2026,-736){\circle*{60}} \put(1876,-736){\makebox(0,0){$R_2$}}
\put(2026,-1411){\circle*{60}} \put(1876,-1411){\makebox(0,0){$R_6$}}
\put(3376,-736){\circle*{60}} \put(3526,-736){\makebox(0,0){$R_1$}}
\put(3376,-1411){\circle*{60}} \put(3526,-1411){\makebox(0,0){$R_3$}}
\put(2701,-1411){\circle*{60}} \put(2851,-1411){\makebox(0,0){$R_{12}$}}
{\color[rgb]{0,0,0}\put(2701,-61){\line( 1,-2){675}}
\put(3376,-1411){\line( 0, 1){675}}
\put(3376,-736){\line(-1, 0){1350}}
\put(2026,-736){\line( 2,-1){1350}}
}%
\end{picture}%
\vspace{1\baselineskip}

Let $G={}^3D_4(q)$. From the compact form for $GK({}^3D_4(q))$ we immediately obtain that
$\Theta(G)=\{\{r_3,r_6,r_{12}\}\mid r_i\in R_i\}$ and $\Theta'(G)=\varnothing$ if $q\not= 2$. For
$q=2$ the result follows from the compact form for the prime graph $GK({}^3D_4(q))$, and the fact
that $R_6=\varnothing$. %\vspace{1\baselineskip}

Let $G={}^2B_2(q)$. In this case primes $s_i\in S_i$ and $s_j\in S_j$ are adjacent if and only if $i=j$, while $p=2$ is non-adjacent to all
vertices, and the proposition follows.

Let $G={}^2G_2(q)$. In this case odd primes $s_i\in S_i$ and $s_j\in S_j$ are adjacent if and only
if $i=j$, while $p=3$ is non-adjacent to all odd primes. Since $2$ is adjacent to $s_1$, $s_2$, and
$p$, we obtain the statement of the proposition in this case.

%\newpage

Let $G={}^2F_4(q)$. If $q>8$, then any set of type $\{s_2,s_3,s_4,s_5,s_6\}$ forms a coclique in
$GK(G)$ by Proposition \ref{adjsuzree}. The same proposition together with
\cite[Proposition~3.3]{VasVd} implies that the set $\{2\}\cup S_1\cup S_2$ forms a clique in
$GK(G)$, any prime $s_3$ is adjacent to $s_1$ and $2$, and $3$ is adjacent to $s_2$ and $s_4$. By
using this information we obtain the proposition in this case. If $G={}^2F_4(8)$ then
$S_2=\pi(9)\setminus\{3\}=\varnothing$, whence every $\theta(G)\in\Theta(G)$ is of type
$\{s_5,s_6\}$, and every $\theta'(G)\in\Theta'(G)$ is a two-element set of type either
$\{s_1,s_4\}$, or  $\{3,s_3\}$, or $\{2,s_4\}$, or $\{s_3,s_4\}$. The group $G={}^2F_4(2)$ is not
simple, and its derived subgroup $T={}^2F_4(2)'$ is the simple Tits group. Using \cite{Atlas}, we
obtain that the prime graph of the Tits group $T$ contains a unique coclique $\rho(T)=\{3,5,13\}$
of maximal size.
\end{proof}

\begin{prop}\label{cocliquea12} If $G\simeq A_{n-1}^\varepsilon(q)$ is an  finite simple groups of Lie type
over a field of characteristic~$p$ and $n\in \{2,3\}$, then $t(G)$, and the sets $\Theta(G)$, $\Theta'(G)$ are listed
in Table~{\em\ref{LinearUnitaryTable}}.
\end{prop}

\begin{proof}
Let $G=A_1(q)$.   Then the compact form for $GK(A_1(q))$ is a coclique with the set of vertices $\{R_1,R_2,p\}$.
Thus $\Theta(G)=\{\{r_1,r_2,p\}\mid r_i\in R_i\}$
and $\Theta'(G)=\varnothing$.
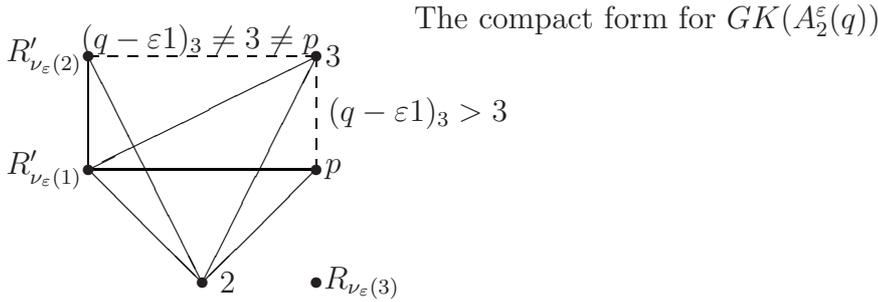
\begin{center}
The compact form for $GK(A_2^\varepsilon(q))$
\end{center}

\setlength{\unitlength}{4144sp}%
\begingroup\makeatletter\ifx\SetFigFontNFSS\undefined%
\gdef\SetFigFontNFSS#1#2#3#4#5{%
  \reset@font\fontsize{#1}{#2pt}%
  \fontfamily{#3}\fontseries{#4}\fontshape{#5}%
  \selectfont}%
\fi\endgroup%
\begin{picture}(1374,1374)(1114,-2323)
\thinlines
{\color[rgb]{0,0,0}\multiput(1126,-961)(117.39130,0.00000){12}{\line( 1, 0){ 58.696}}
}%
{\color[rgb]{0,0,0}\put(1126,-961){\line( 1,-2){675}}
}%
{\color[rgb]{0,0,0}\put(1801,-2311){\line( 1, 2){675}}
}%
{\color[rgb]{0,0,0}\put(2476,-1636){\line(-1,-1){675}}
}%
{\color[rgb]{0,0,0}\put(1126,-1636){\line( 1,-1){675}}
}%
{\color[rgb]{0,0,0}\put(1126,-1636){\line( 1, 0){1350}}
}%
{\color[rgb]{0,0,0}\put(1126,-961){\line( 0,-1){675}}
}%
{\color[rgb]{0,0,0}\multiput(2476,-961)(0.00000,-122.72727){6}{\line( 0,-1){ 61.364}}
}%
{\color[rgb]{0,0,0}\put(2476,-961){\line(-2,-1){1350}}
}%
\put(1126,-961){\circle*{60}} \put(866,-961){\makebox(0,0){$R_{\nu_\varepsilon(2)}'$}}
\put(1786,-861){\makebox(0,0){$(q-\varepsilon1)_3\not=3\not=p$}}
\put(3076,-1300){\makebox(0,0){$(q-\varepsilon1)_3>3$}}
\put(1801,-2311){\circle*{60}} \put(1951,-2311){\makebox(0,0){$2$}}
\put(2476,-961){\circle*{60}} \put(2576,-961){\makebox(0,0){$3$}}
\put(2476,-1636){\circle*{60}} \put(2576,-1636){\makebox(0,0){$p$}}
\put(1126,-1636){\circle*{60}} \put(866,-1636){\makebox(0,0){$R_{\nu_\varepsilon(1)}'$}}
\put(2476,-2311){\circle*{60}} \put(2746,-2311){\makebox(0,0){$R_{\nu_\varepsilon(3)}$}}
\end{picture}%
\vspace{1\baselineskip}

Let $G=A_2^\varepsilon(q)$. Set $R_{\nu_\varepsilon(1)}'= R_{\nu_\varepsilon(1)}\setminus\{2,3\}$,
and $R_{\nu_\varepsilon(2)}'= R_{\nu_\varepsilon(2)}\setminus\{2,3\}$.  Assume first that
$(q-\varepsilon1)_3>3$. Then the set $\{2,3,p,R_{\nu_\varepsilon(1)}'\}$ is a clique in the compact
form for $GK(A_2^\varepsilon(q))$, while $R_{\nu_\varepsilon(2)}'$ and $R_{\nu_\varepsilon(3)}$ are
non-adjacent. If $R_{\nu_\varepsilon(2)}\not=\{2\}$ (i.~e., $q+\varepsilon1\not=2^k$ and
$R_{\nu_\varepsilon(2)}'\not=\varnothing$), then $p$ is the only vertex from the clique
$\{2,3,p,R_{\nu_\varepsilon(1)}'\}$, which is non-adjacent to both $R_{\nu_\varepsilon(2)}'$ and
$R_{\nu_\varepsilon(3)}$. Hence
$\Theta(G)=\{\{p,r_{\nu_\varepsilon(2)}\not=2,r_{\nu_\varepsilon(3)}\}\mid r_i\in R_i\}$ and
$\Theta'(G)=\varnothing$. If $R_{\nu_\varepsilon(2)}=\{2\}$ (i.~e., $q+\varepsilon1=2^k$ and
$R_{\nu_\varepsilon(2)}'=\varnothing$), then $\Theta(G)=\{\{r_{\nu_\varepsilon(3)}\}\mid
r_{\nu_\varepsilon(3)}\in R_{\nu_\varepsilon(3)}\}$ and
$\Theta'(G)=\{\{2\},\{p\},\{r_{\nu_\varepsilon(1)}\}\mid r_{\nu_\varepsilon(1)}\in
R_{\nu_\varepsilon(1)} \}$.

Now assume that $(q-\varepsilon1)_3=3$. Then the set $\{2,p,R_{\nu_\varepsilon(1)}'\}$ is a clique
in the compact form for $GK(A_2^\varepsilon(q))$, while $3$, $R_{\nu_\varepsilon(2)}'$, and
$R_{\nu_\varepsilon(3)}$ are pairwise non-adjacent. Since $p$ is the only vertex from the clique
$\{2,p,R_{\nu_\varepsilon(1)}'\}$, which is non-adjacent to  $3$, $R_{\nu_\varepsilon(2)}'$, and
$R_{\nu_\varepsilon(3)}$, we obtain that
$\Theta(G)=\{\{3,p,r_{\nu_\varepsilon(2)}\not=2,r_{\nu_\varepsilon(3)}\}\mid r_i\in R_i\}$ if
$R_{\nu_\varepsilon(2)}\not=\{2\}$, and $\Theta(G)=\{\{3,p,r_{\nu_\varepsilon(3)}\}\mid
r_{\nu_\varepsilon(3)}\in R_{\nu_\varepsilon(3)}\}$ if $R_{\nu_\varepsilon(2)}=\{2\}$. In both
cases $\Theta'(G)=\varnothing$.

Assume at the end that $(q-\varepsilon1)_3=1$, i.~e., either $(q+\varepsilon1)_3>1$ and  $3\in
R_{\nu_\varepsilon(2)}\not=\{2\}$, or $p=3$. As above we have that the set
$\{2,p,R_{\nu_\varepsilon(1)}'\}$ is a clique in the compact form for $GK(A_2^\varepsilon(q))$,
while $R_{\nu_\varepsilon(2)}'$ and $R_{\nu_\varepsilon(3)}$ are pairwise non-adjacent. Since $p$
is the only vertex from the clique $\{2,p,R_{\nu_\varepsilon(1)}'\}$, which is non-adjacent to
$R_{\nu_\varepsilon(2)}'$ and $R_{\nu_\varepsilon(3)}$, and since either $3\in
R_{\nu_\varepsilon(2)}$ or $p=3$, we obtain that
$\Theta(G)=\{\{p,r_{\nu_\varepsilon(2)}\not=2,r_{\nu_\varepsilon(3)}\}\mid r_i\in
R_i\setminus\{2\}\}$ and $\Theta'(G)=\varnothing$ if $R_{\nu_\varepsilon(2)}\not=\{2\}$,  and
$\Theta(G)=\{\{r_{\nu_\varepsilon(3)}\}\mid r_{\nu_\varepsilon(3)}\in R_{\nu_\varepsilon(3)}\}$ and
$\Theta'(G)=\{\{p\},\{r_{\nu_\varepsilon(1)}\},\{2=r_{\nu_\varepsilon(2)}\}\mid
r_{\nu_\varepsilon(1)}\in R_{\nu_\varepsilon(1)}\}$ if $R_{\nu_\varepsilon(2)}=\{2\}$.
\end{proof}

Below we give Tables~\ref{LinearUnitaryTable}, \ref{ClassicTable}, \ref{ExceptTable}. These tables
are organized in the following way. Column~1 represents a group of Lie type $G$ with the base field
of order $q$ and characteristic $p$, Column~2 contains conditions on $G$, and Column~3 contains
value of $t(G)$. In Columns~4 and~5 we list the elements of $\Theta(G)$ and $\Theta'(G)$, that is
sets $\theta(G)\in\Theta(G)$ and $\theta'(G)\in\Theta'(G)$, and omit the braces for one-element
sets. In particular, the item $\{p,3,r_2\not=2,r_3\}$ in Column~4 means
$\Theta(G)=\{\{p,3,r_2,r_3\}\mid r_2\in R_2\setminus\{2\},r_3\in R_3\}$ and the item $p,r_4$ in
Column~5 means $\Theta'(G)=\{\{p\},\{r_4\}\mid r_4\in R_4\}$.

\newpage

\begin{tab}\label{LinearUnitaryTable}{\bfseries Cocliques for finite simple linear and unitary groups}
\smallskip

{\small    \noindent\begin{tabular}{|c|l|c|c|c|}
  \hline
  $G$ & Conditions & $t(G)$ & $\Theta(G)$ & $\Theta'(G)$\\
  \hline
$A_1(q)$&$q>3$&$3$&$\{p,r_1,r_2\}$&$\varnothing$\\ \hline
$A_2(q)$&$(q-1)_3=3$, $q+1\not=2^k$& $4$&$\{p,3,r_2\not=2,r_3\}$&$\varnothing$ \\
& $(q-1)_3=3$, $q+1=2^k$& $3$&$\{3,p,r_3\}$&$\varnothing$ \\
&$(q-1)_3\not=3$, $q+1\not=2^k$& $3$&$\{p,r_2\not=2,r_3\}$&$\varnothing$ \\
&$(q-1)_3\not=3$, $q+1=2^k$& $2$&$r_3$&$p,r_1,2=r_2$ \\ \hline
$A_3(q)$&$(q-1)_2\not=4$&$3$&$\{p,r_3,r_4\}$&$\varnothing$\\
&$(q-1)_2=4$&$3$&$\{r_3,r_4\}$&$p,2$\\ \hline
$A_4(q)$& $(q-1)_5\neq5$ &$3$&$\{r_4,r_5\}$&$p,r_3$\\
& $(q-1)_5=5$ &$3$&$\{r_4,r_5\}$&$5,p,r_3$\\
\hline
 $A_5(q)$&$q=2$&$3$&$\{r_3,r_4,r_5\}$&$\varnothing$\\
 &$q>2$ and $(q-1)_3\neq3$ &$3$&$r_5$&$\{p,r_6\}$,$\{r_3,r_4\}$, \\
 & & & & $\{r_4,r_6\}$\\
 &$(q-1)_3=3$ &$3$&$r_5$&$\{p,r_6\}$,$\{r_3,r_4\}$, \\
 & & & & $\{r_4,r_6\}$, $\{3,r_6\}$\\ \hline
$A_{n-1}(q),$ &$n$ is  odd and $q\not=2$&$[\frac{n+1}{2}]$&$\{r_i\mid \frac{n}{2}< i\le n\}$&$\varnothing$\\
$n\ge7$& for $7\le n\le 11$&&&\\
& $n$ is even and
 $q\not=2$&$[\frac{n+1}{2}]$&$\{ r_{i}\mid \frac{n}{2}< i<
n\}$&$r_{\frac{n}{2}}$,$r_n$\\
& for $8\le n\le12$&&&\\
&$n=7$, $q=2$&$3$&$\{r_5,r_7\}$&$r_3,r_4$\\
&$n=8$, $q=2$&$3$&$r_7$& $\{p,r_8\}$,$\{r_5,r_8\}$, \\
&&&&$\{r_3,r_8\}$,$\{r_4,r_5\}$\\
&$n=9$, $q=2$&$4$&$\{r_5,r_7,r_8,r_9\}$&$\varnothing$\\
&$n=10$, $q=2$&$4$&$\{r_7,r_9\}$&$\{r_4,r_{10}\}$,$\{r_8,r_{10}\}$\\
&&&&$\{r_5,r_{8}\}$\\
&$n=11$, $q=2$&$5$&$\{r_7,r_8,r_9,r_{11}\}$&$r_5,r_{10}$ \\
&$n=12$, $q=2$&$6$&$\{r_7,r_8,r_9,r_{10},r_{11},r_{12}\}$&$\varnothing$\\ \hline
${}^2A_2(q)$,&$(q+1)_3=3$, $q-1\not=2^k$&$4$&$\{p,3,r_{1}\not=2,r_{6}\}$&$\varnothing$\\
$q>2$&$(q+1)_3=3$, $q-1=2^k$&$3$&$\{3,p,r_{6}\}$&$\varnothing$\\
&$(q+1)_3\not=3$, $q-1\not=2^k$&$3$&$\{p,r_{1}\not=2,r_{6}\}$&$\varnothing$\\
&$(q+1)_3\not=3$, $q-1=2^k>2$&$2$&$r_{6}$&$p, r_{2},2=r_{1}$\\
&$q=3$&$2$&$r_{6}$&$p, r_{2}=2$\\ \hline
${}^2A_3(q)$&$(q+1)_2\not=4$ and $q\neq2$ &$3$&$\{p,r_{6},r_{4}\}$&$\varnothing$\\
&$(q+1)_2=4$&$3$&$\{r_{6},r_{4}\}$&$p,2$\\
& $q=2$ & $2$ & $r_4$ & $p,r_2$\\ \hline
 ${}^2A_4(q)$& $q=2$ &$3$&$\{p,r_{4},r_{10}\}$&$\varnothing$\\
 & $q>2$ and $(q+1)_5\neq5$ &$3$&$\{r_{4},r_{10}\}$&$p,r_{6}$\\
 & $(q+1)_5=5$ &$3$&$\{r_{4},r_{10}\}$&$5,p,r_{6}$\\
\hline
 ${}^2A_5(q)$&$q=2$&$3$&$\{r_{10},r_{3}\}$&$3,p,r_{4}$\\
 & $(q+1)_3\neq3$&$3$&$r_{10}$&$\{p,r_{3}\}$,$\{r_{6},r_{4}\}$, \\
 & & & & $\{r_{4},r_{3}\}$\\
 & $q>2$ and $(q+1)_3=3$&$3$&$r_{10}$&$\{p,r_{3}\}$,$\{r_{6},r_{4}\}$, \\
 & & & & $\{r_{4},r_{3}\}$, $\{3,r_{3}\}$\\
 \hline
 ${}^2A_{n-1}(q)$,&$n$ is odd&$[\frac{n+1}{2}]$&$\{r_i\mid \frac{n}{2}< \nu(i)\le n\}$&$\varnothing$\\
$n\ge7$&$n$ is even &$[\frac{n+1}{2}]$&$\{r_i\mid \frac{n}{2}< \nu(i)< n\}$&$r_{\nu(\frac{n}{2})},r_{\nu(n)}$\\
\hline
\end{tabular}}
\end{tab}
\newpage

\begin{tab}\label{ClassicTable}{\bfseries Cocliques for finite simple symplectic and orthogonal groups}
\smallskip

{\small    \noindent\begin{tabular}{|c|l|c|c|c|}
  \hline
  $G$ & Conditions & $t(G)$ & $\Theta(G)$ & $\Theta'(G)$\\
  \hline
 $B_n(q)$ or & $n=2$, $q=3$ & 2 & $r_4$ & $p,r_2$\\
 $C_n(q)$ & $n=2$, $q>3$ & 2 & $r_4$ & $p,r_1,r_2$\\
 & $n=3$ and $q=2$ & $2$ & $r_3$ & $p,r_2,r_4$\\
 & $n=3$ and $q>2$ & $3$ & $\{r_3,r_6\}$ & $p,r_4$\\
 & $n=4$ and $q=2$ & $3$ & $\{r_3,r_4,r_8\}$ & $\varnothing$\\
 & $n=5$ and $q=2$ & $4$ & $\{r_5,r_8,r_{10}\}$ & $r_3,r_4$\\
 & $n=6$ and $q=2$ & $5$ & $\{r_3,r_5,r_8,r_{10},r_{12}\}$ & $\varnothing$\\
 & $n=7$ and $q=2$ & $6$ & $\{r_5,r_7,r_{10},r_{12},r_{14}\}$ & $r_3,r_8$\\
 & $n>3$, $n\equiv{0,1}(\mod 4)$ and & $\left[\frac{3n+5}{4}\right]$ & $\{r_i\mid
 \frac{n}{2}\leqslant\eta(i)\leqslant n\}$ &
 $\varnothing$\\
 & $(n,q)\neq(4,2),(5,2)$ &  & &\\
 & $n>3$, $n\equiv{2}(\mod 4)$ and & $\left[\frac{3n+5}{4}\right]$ & $\{r_i\mid \frac{n}{2}<\eta(i)\leqslant n\}$ & $r_{n/2},r_n$\\
 & $(n,q)\neq(6,2)$ &  & &\\
 & $n>3$, $n\equiv{3}(\mod 4)$ and & $\left[\frac{3n+5}{4}\right]$ & $\{r_i\mid \frac{n+1}{2}<\eta(i)\leqslant n\}$ &
 $r_{(n-1)/2},r_{n-1},$\\
 & $(n,q)\neq(7,2)$ &  & & $r_{n+1}$\\
 \hline
 $D_n(q)$ & $n=4$, $q=2$ & 2 & $r_3$ & $p,r_2,r_4$\\
 & $n=4$ and $q>2$ & $3$ & $\{r_3,r_6\}$ & $p,r_4$\\
 & $n=5$ and $q=2$ & $4$ & $\{r_3,r_4,r_5,r_8\}$ & $\varnothing$\\
 & $n=6$ and $q=2$ & $4$ & $\{r_3,r_5,r_8,r_{10}\}$ & $\varnothing$\\
 & $n>4$, $n\equiv{0}(\mod 4)$ & $\left[\frac{3n+1}{4}\right]$ & $\{r_i\mid \frac{n}{2}\leqslant\eta(i)\leqslant n,$ & $\varnothing$\\
 & & & $i\neq2n\}$ & \\
 & $n>4$, $n\equiv{1}(\mod 4)$ and & $\left[\frac{3n+1}{4}\right]$ & $\{r_i\mid  \frac{n}{2}<\eta(i)\leqslant n,$ & $r_{n-1},r_{n+1}$\\
 & $(n,q)\neq(5,2)$ &  & $i\neq2n,n+1\}$ &\\
 & $n>4$, $n\equiv{2}(\mod 4)$ and & $\left[\frac{3n+1}{4}\right]$ & $\{r_i\mid \frac{n}{2}<\eta(i)\leqslant n,$ & $r_{n/2},r_{n}$\\
 & $(n,q)\neq(6,2)$ &  & $i\neq2n\}$ &\\
 & $n>4$, $n\equiv{3}(\mod 4)$ & $\frac{3n+3}{4}$ & $\{r_i\mid \frac{n-1}{2}\leqslant\eta(i)\leqslant n,$ & $\varnothing$\\
 & & & $i\neq2n,n-1\}$ &\\
 \hline
 ${}^2D_n(q)$ & $n=4$, $q=2$ & 3 & $\{r_3,r_8\}$ & $p,r_4$\\
 & $n=4$ and $q>2$ & $4$ & $\{r_3,r_6,r_8\}$ & $p,r_4$\\
 & $n=5$ and $q=2$ & $3$ & $\{r_8,r_{10}\}$ & $p,r_3,r_4$\\
 & $n=6$ and $q=2$ & $5$ & $\{r_5,r_8,r_{10},r_{12}\}$ & $r_3,r_4$\\
 & $n=7$ and $q=2$ & $5$ & $\{r_5,r_{10},r_{12},r_{14}\}$ & $r_3,r_8$\\
 & $n>4$, $n\equiv{0}(\mod 4)$ and & $\left[\frac{3n+4}{4}\right]$ & $\{r_i\mid \frac{n}{2}\leqslant\eta(i)\leqslant n\}$ & $\varnothing$\\
 & $n>4$, $n\equiv{1}(\mod 4)$ and & $\left[\frac{3n+4}{4}\right]$ & $\{r_i\mid \frac{n}{2}<\eta(i)\leqslant n,$
 & $r_{(n+1)/2},r_{n-1}$\\
 & $(n,q)\neq(5,2)$ & & $i\neq n,\frac{n+1}{2}\}$ & \\
 & $n>4$, $n\equiv{2}(\mod 4)$ and & $\left[\frac{3n+4}{4}\right]$ & $\{r_i\mid \frac{n}{2}<\eta(i)\leqslant n\}$ &
 $r_{n/2},r_{n-2},r_n$\\
 & $(n,q)\neq(6,2)$ &  &  & \\
 & $n>4$, $n\equiv{3}(\mod 4)$ and & $\left[\frac{3n+4}{4}\right]$ & $\{r_i\mid \frac{n-1}{2}\leqslant\eta(i)\leqslant n,$ & $\varnothing$\\
 & $(n,q)\neq(7,2)$ &  &$i\neq n,\frac{n-1}{2}\}$ &\\
 \hline
\end{tabular}}
\end{tab}

\newpage

\begin{tab}\label{ExceptTable}{\bfseries Cocliques for finite simple exceptional groups}
\smallskip

{\small    \noindent\begin{tabular}{|c|l|c|c|c|}
  \hline
  $G$ & Conditions & $t(G)$ & $\Theta(G)$ & $\Theta'(G)$\\
  \hline
  $G_2(q)$ & $q=3,4$ &3&  $\{r_3,r_6\}$ & $p,r_2$ \\
  & $q=8$ &3&  $\{r_3,r_6\}$ & $p,r_1$ \\
  & $q=3^m>3$ &3&  $\{r_3,r_6\}$ & $p,r_1,r_2$ \\
  & $q\equiv{1}(\mod 3)$ and $q\neq4$ & 3 & $\{r_3,r_6\}$ & $p,r_2,r_1\neq3$\\
  & $q\equiv{2}(\mod 3)$ and $q\neq8$ & 3 & $\{r_3,r_6\}$ & $p,r_1,r_2\neq3$\\ \hline
  $F_4(q)$&$q=2$&4&$\{r_3,r_4,r_8,r_{12}\}$ & $\varnothing$ \\
  & $q>2$&5&$\{r_3,r_4,r_6,r_8,r_{12}\}$ & $\varnothing$ \\ \hline
  $E_6(q)$&$q=2$&5&$\{r_4,r_5,r_8,r_9\}$ & $r_3,r_{12}$\\
  &$q>2$&5&$\{r_5,r_8,r_9\}$ & $\{r_3,r_4\}$, $\{r_4,r_{12}\}$, \\
  & & & & $\{r_6,r_{12}\}$ \\ \hline
  ${}^2E_6(q)$ & $q=2$ & 5 & $\{r_8,r_{10},r_{12},r_{18}\}$ & $r_3,r_4$\\
  & $q>2$ & 5 & $\{r_8,r_{10},r_{18}\}$ & $\{r_3,r_{12}\}$, $\{r_4,r_6\}$, \\
  & & & & $\{r_4,r_{12}\}$ \\
  \hline
  $E_7(q)$&&8&$\{r_5,r_7,r_9,r_{10},$ & $r_4,r_8$\\
  &&&$r_{12},r_{14},r_{18}\}$&\\  \hline
  $E_8(q)$&&12&$\{r_5,r_7,r_8,r_9,r_{10},r_{12},$ & $\varnothing$\\
  &&&$r_{14},r_{15},r_{18},r_{20},r_{24},r_{30}\}$&\\
  \hline
 ${}^3D_4(q)$& $q=2$ & 2 & $r_{12}$ & $p,r_2,r_3$\\
 & $q>2$ & 3 & $\{r_3,r_6,r_{12}\}$ & $\varnothing$\\
 \hline
 ${^2B_2(2^{2n+1})}$&$n\ge1$&4&$\{2,s_1,s_2,s_3\}$ & $\varnothing$\\
  \hline
 ${^2G_2(3^{2n+1})}$&$n\ge1$&5&$\{3,s_1,s_2,s_3,s_4\}$ & $\varnothing$\\
 \hline
 ${}^2F_4(2^{2n+1})$&$n\ge2$,&$5$&$\{s_2,s_3,s_4,s_5,s_6\}$  & $\varnothing$\\
 \hline
 ${}^2F_4(8)$ & &  $4$&$\{s_5,s_6\}$ & $\{3,s_3\}$, $\{s_1,s_4\}$,\\
 &&&&$\{2,s_4\}$, $\{s_3,s_4\}$\\\hline

 ${}^2F_4(2)'$ & &  $3$&$\{3,5,13\}$ & $\varnothing$\\ \hline

\end{tabular}}
\end{tab}

\newpage

\section{Appendix}\label{appendix}

In this section we give a list of corrections for~\cite{VasVd} which we obtain in the present paper.

Items (4), (5), (9) of Lemma 1.3 should be substituted by items (1), (2), (3) of Lemma
\ref{toriofexcptgrps} of the present paper respectively.

Lemma 1.4 should be substituted by Lemma \ref{Zsigmondy Theorem}.

Lemma 1.5 should be substituted be Lemma \ref{SuzReeDivisors}.

Proposition 2.3 should be substituted by Proposition \ref{adjbn}.

Proposition 2.4 should be substituted by Proposition \ref{adjdn}.

Proposition 2.5 should be substituted by Proposition \ref{adjexcept}.

In Tables 4 and 8  the following corrections are necessary.

The lines

\begin{longtable}{|c|c|c|c|}
$A_{n-1}(q)$&$n=3$, $(q-1)_3=3$, and $q+1\not=2^k$&$4$&$\{p,3,r_2,r_3\}$\\
&$n=3$, $(q-1)_3\not=3$, and $q+1\not=2^k$&$3$&$\{p,r_2,r_3\}$\\
\end{longtable}

should be substituted by the lines

\begin{longtable}{|c|c|c|c|}
$A_{n-1}(q)$&$n=3$, $(q-1)_3=3$, and $q+1\not=2^k$&$4$&$\{p,3,r_2\not=2,r_3\}$\\
&$n=3$, $(q-1)_3\not=3$, and $q+1\not=2^k$&$3$&$\{p,r_2\not=2,r_3\}$\\
\end{longtable}

The lines

\begin{longtable}{|c|c|c|c|}
${}^2A_{n-1}(q)$&$n=3$, $(q+1)_3=3$, and $q-1\not=2^k$&$4$&$\{p,3,r_1,r_6\}$\\
&$n=3$, $(q+1)_3\not=3$, and $q-1\not=2^k$&$3$&$\{p,r_1,r_6\}$\\
\end{longtable}

should be substituted by the lines

\begin{longtable}{|c|c|c|c|}
${}^2A_{n-1}(q)$&$n=3$, $(q+1)_3=3$, and $q-1\not=2^k$&$4$&$\{p,3,r_1\not=2,r_6\}$\\
&$n=3$, $(q+1)_3\not=3$, and $q-1\not=2^k$&$3$&$\{p,r_1\not=2,r_6\}$\\
\end{longtable}

In Table 4 in the penultimate line corresponding to  $D_n(q)$ instead of
$n\equiv 1\pmod 1$, $n>4$ there should be  $n\equiv 1\pmod 2$, $n>4$.

In Table 8 the following corrections are necessary.

The line

\begin{longtable}{|c|c|c|c|}
$D_n(q)$&$n\ge 4$, $(n,q)\not=(4,2),(5,2),(6,2)$&$\left[\frac{3n+1}{4}\right]$&$\{r_{2i}\mid \left[\frac{n+1}{2}\right]\le
i< n\}\cup$\\
&&&$\cup\{r_i\mid \left[\frac{n}{2}\right]<i\le
n,$\\
&&&$i\equiv 1\pmod2\}$
\end{longtable}

should be substituted by

\begin{longtable}{|c|c|c|c|}
$D_n(q)$&$n\ge 4$, $n\not\equiv3\pmod4$,&$\left[\frac{3n+1}{4}\right]$&$\{r_{2i}\mid
\left[\frac{n+1}{2}\right]\le
i< n\}\cup$\\
&$(n,q)\not=(4,2),(5,2),(6,2)$&&$\cup\{r_i\mid \left[\frac{n}{2}\right]<i\le
n,$\\
&&&$i\equiv 1\pmod2\}$\\
&$n\equiv3\pmod4$&$\frac{3n+3}{4}$&$\{r_{2i}\mid \left[\frac{n+1}{2}\right]\le
i< n\}\cup$\\
&&&$\cup\{r_i\mid \left[\frac{n}{2}\right]\le i\le
n,$\\
\end{longtable}

In Table 8 the line

\begin{longtable}{|c|c|c|c|}
${}^2D_n(q)$&$n\ge 4$, $n\not\equiv1\pmod
4,$&$\left[\frac{3n+4}{4}\right]$&$\{r_{2i}\mid \left[\frac{n}{2}\right]\le
i\le n\}\cup$\\
&$(n,q)\not=(4,2),(6,2),(7,2)$&&$\cup\{r_i\mid \left[\frac{n}{2}\right]<i\le
n,$\\
&&&$i\equiv 1\pmod2\}$
\end{longtable}

should be substituted by the line

\begin{longtable}{|c|c|c|c|}
${}^2D_n(q)$&$n\ge 4$, $n\not\equiv1\pmod
4,$&$\left[\frac{3n+4}{4}\right]$&$\{r_{2i}\mid \left[\frac{n}{2}\right]\le
i\le n\}\cup$\\
&$(n,q)\not=(4,2),(6,2),(7,2)$&&$\cup\{r_i\mid \left[\frac{n}{2}\right]<i<
n,$\\
&&&$i\equiv 1\pmod2\}$
\end{longtable}

In Table 9 the following corrections are necessary.

The line

\begin{longtable}{|c|c|c|c|}\hline
$E_6(q)$&$q=2$&5&$\{5,12,17,19,31\}$\\
&$q>2$&$6$&$\{r_4,r_5,r_6,r_8,r_9,r_{12}\}$\\ \hline
\end{longtable}

should be substituted by the line %\newpage

\begin{longtable}{|c|c|c|c|}\hline
$E_6(q)$&none&5&$\{r_4,r_5,r_8,r_9,r_{12}\}$\\ \hline
\end{longtable}

The line

\begin{longtable}{|c|c|c|c|}\hline
$E_7(q)$&none&7&$\{r_7,r_8,r_9,r_{10},r_{12},r_{14},r_{18}\}$\\ \hline
\end{longtable}

should be substituted by the line

\begin{longtable}{|c|c|c|c|}\hline
$E_7(q)$&none&8&$\{r_5,r_7,r_8,r_9,r_{10},r_{12},r_{14},r_{18}\}$\\ \hline
\end{longtable}

The line

\begin{longtable}{|c|c|c|c|}\hline
$E_8(q)$&none&11&$\{r_7,r_8,r_9,r_{10},r_{12},r_{14},r_{15},r_{18},r_{20},r_{24},
r_{30}\} $\\ \hline
\end{longtable}

should be substituted by the line

\begin{longtable}{|c|c|c|c|}\hline
$E_8(q)$&none&12&$\{r_5,r_7,r_8,r_9,r_{10},r_{12},r_{14},r_{15},r_{18},r_{20},r_{
24 } ,
r_{30}\} $\\ \hline
\end{longtable}

A revised variant of \cite{VasVd} can be found in

http://arxiv.org/abs/math/0506294.

\end{document}